\documentclass[12pt]{article}
\usepackage{amsmath}
\usepackage{amssymb}
\usepackage{harpoon}
\usepackage{latexsym,color}
\usepackage[textures,dvipdfm]{hyperref}
\usepackage{ntheorem}
\usepackage{mathrsfs}
\usepackage[numbers,sort&compress]{natbib}

  \allowdisplaybreaks[4]
\title{ \Large On the convergence to equilibrium for the spatially homogeneous Boltzmann equation for Fermi-Dirac particles}
\author{ Bocheng Liu\footnote{ Department of Mathematical Sciences, Tsinghua
University, Beijing 100084, P.R.China; \,\, e-mail address: liu-bc16@mails.tsinghua.edu.cn }\,\,
and  Xuguang Lu\footnote{Department of Mathematical Sciences,
Tsinghua University, Beijing 100084, P.R.China;  e-mail address: xglu@mail.tsinghua.edu.cn
} }

\date{}  

\newcommand{\td}{\tilde}
\newcommand{\ld}{\lambda}

\newcommand{\p}{\partial}
\newcommand{\vp}{\varphi}
\newcommand{\vep}{\varepsilon}
\newcommand{\og}{\omega}

\newcommand{\sg}{\sigma}

\newcommand{\gm}{\gamma}
\newcommand{\Gm}{\Gamma}
\newcommand{\dt}{\delta}
\newcommand{\Dt}{\Delta}
\newcommand{\fr}{\frac}
\newcommand{\wt}{\widetilde}
\newcommand{\wh}{\widehat}

\newcommand{\bR}{{\mathbb R}^3}
\newcommand{\bS}{{\mathbb S}^2 }
\newcommand{\bSS}{{\bS}\times{\bS}}

\newcommand{\bRR}{{\bR}\times{\bR}}

\newcommand{\bRS}{{\bR}\times {\mathbb S}^2 }
\newcommand{\bRRS}{{\bRR}\times{\mathbb S}^2 }

\newcommand{\la}{\langle}
\newcommand{\ra}{\rangle}
\newcommand{\mR}{{\mathbb R}}
\newcommand{\mN}{{\mathbb N}}

\newcommand{\mS}{{\mathbb S}}
\newcommand{\be}{\begin{myequation}}
\newcommand{\ee}{\end{myequation}}
\newcommand{\bes}{\begin{myeqnarray}}
\newcommand{\ees}{\end{myeqnarray}}
\newcommand{\beas}{\begin{eqnarray*}}
\newcommand{\eeas}{\end{eqnarray*}}
\newcommand{\lb}{\label}
\newcounter{thm}
\theoremseparator{.}
\newtheorem{theorem}{Theorem}[section]
\newtheorem{proposition}[theorem]{Proposition}

\newtheorem{lemma}[theorem]{Lemma}
\newtheorem{remark}[theorem]{Remark}

\newcounter{myequation}[section]

\newenvironment{proof}{{\bf Proof.}}{$\hfill\Box$}
\newenvironment{myequation}{\stepcounter{myequation}\begin{equation}}{\end{equation}}
\newenvironment{myeqnarray}{\stepcounter{myequation}\begin{eqnarray}}{\end{eqnarray}}
\newcommand{\dnumber}{\stepcounter{myequation}}

\begin{document}
\maketitle

\vskip 0.1in \baselineskip 18pt
\begin{abstract}
	
	In this paper we prove the strong and time-averaged strong convergence to equilibrium for solutions
(with general initial data) of the spatially homogeneous Boltzmann equation for Fermi-Dirac particles. The assumption on the collision kernel includes the Coulomb potential with a weaker angular cutoff. The proof is based on moment estimates, entropy dissipation inequalities, regularity of the collision gain operator, and a new observation that many collision kernels are larger than or equal to some completely positive kernels, which enables us to avoid dealing with the convergence problem of the cubic collision integrals.


\end{abstract}

\begin{center}\section { Introduction}\end{center}
The quantum Boltzmann equations, in which the collision integrals are defined in terms of Bose-Einstein statistics and Fermi-Dirac statistics,
 describe mesoscopic behaviour of dilute gases of bosons and fermions respectively.
These equations were first derived by Nordheim \cite{Nordheim} and Uehling $\&$ Uhlenbeck \cite{Uehling and Uhlenbeck}
and then developed by \cite{weak-coupling},\cite{Chapman and Cowling},\cite{ESY},\cite{EMV0},\cite{LS}.
In compared with the classical Boltzmann equation, the quantum Boltzmann equations have highly non-linear structures, and this leads to lacking deeper and quantitative results  even for
 the spatially homogeneous solutions so far.
In this paper we consider the spatially homogeneous
 Boltzmann equation for fermions and prove the long time strong convergence to equilibrium.
(Research on this area has
 been no result for the spatially inhomogeneous solutions with general initial (-boundary) data,  except the results on existence of solutions, see e.g. \cite{Do},\cite{Lions2},\cite{Alex},\cite{Lu2008}). After normalized some physical constants, the equation under consideration is written
$$
\frac{\p}{\p t}f(t,v)=Q(f)(t,v),\quad (t,v)\in(0,\infty)\times{\bR}\eqno({\rm BFD})$$
where
\bes&& Q(f)(t,v)=\int_{\bRS}B(v-v_*,\sg)\Pi_{{\rm F}}(f)(t,v,v_*,\sg)
{\rm d}\sg {\rm d}v_*,\lb{Q}\\
&&
\Pi_{{\rm F}}(f)=f'f'_*(1-f)(1-f_*)-ff_*(1-f')(1-f_*'),\dnumber \lb{1.2}\ees
${\rm d}\sg$ is the area element of the sphere measure on ${\mS}^2$,
$f=f(t,v), f_*=f(t,v_*), f'=f(t,v'),f_*'=f(t,v_*')$ , and $v,v_*$ and $v',v_*'$ are velocities of two particles before and after
their collision, and their relation can be expressed as the $\sg$-representation:
\be v'=\fr{v+v_*}{2}+\fr{|v- v_*|\sg}{2},\quad v'_*=\fr{v+v_*}{2}-\fr{|v- v_*|\sg}{2},
\qquad \sg\in{\mathbb S}^2
\lb{1.1}\ee
which conserves the momentum and kinetic energy
\be v'+v_*'=v+v_*,\quad |v'|^2+|v_*'|^2=|v|^2+|v_*|^2.\lb{1.4}\ee
The solution $f(t,v)$ of Eq.(BFD) is a particle number density function (or the velocity distribution function) of a spatially homogeneous one-species Fermi-Dirac particle system. The quantum effect in this model is
mainly reflected in the shrink factor $0\le 1-f(t,v)\le 1$ (due to the Pauli exclusion principle) which is a part of definition of solutions of Eq.(BFD).
It is not difficult to show that if initially $0\le f|_{t=0}\le 1$ on ${\bR}$, then $0\le f \le 1$ on $[0,\infty)\times {\bR}$.
Although solutions $f$  are $L^{\infty}$-bounded, the non-linear structure $f'f'_*(1-f)(1-f_*)-ff_*(1-f')(1-f'_*)$ in Eq.(BFD) contributes a lot of difficulties in proving
long time strong convergence to equilibrium. In fact it is different from the classical Boltzmann equation
and the Landau-Fermi-Dirac equation (see e.g. \cite{Villani2}, \cite{ABDL} and references therein)
that so far to our knowledge there has been only one partial qualitative result for Eq.(BFD) on long time strong convergence to equilibrium for general initial data (\cite{LW2003}); and even for solutions with
(very) special initial data, results are also quite few (see e.g.\cite{OW}).

As usual we assume that the interaction potential of particles is real and spherically symmetric.
In this case,  the collision kernel $B(z,\sg)$ (or the scattering cross section $B(z,\sg)/|z|$) is a nonnegative Borel function of $|z|$ and $\cos\theta=\la z,\sg\ra/|z|$ only:
\be B(z,\sg)=B(|z|, \cos\theta), \quad
(z,\sg)\in {\bRS}.\lb{1.5}\ee
In some literatures of kinetic theory, it is also used the $\og$-representation:
\be v'=v-\la v-v_*,\og\ra\og,\quad v_*'=v_*+\la v-v_*,\og\ra\og,\quad \og\in{\bS}.
\lb{1.8}\ee
The collision kernel $\wt{B}(z,\og)=\wt{B}(|z|, |\la {\bf n},\og\ra|)$ in Eq.(BFD) corresponding to (\ref{1.8})
has the relation with $B(z,\sg)=B(|z|, \la {\bf n},\sg\ra)$
as follows: $ {\bf n}=z/|z|, \sg={\bf n}-2\la {\bf n},\og\ra\og,$
 $$\wt{B}(|z|, |\la {\bf n},\og\ra|)=2|\la {\bf n},\og\ra|B(|z|,1-2\la {\bf n},\og\ra^2) $$
and it holds the identity:
 \beas&& \int_{{\bS}} B(v-v_*,\sg)\psi(v',v_*')\big|_{\sg-{\rm rep.}}{\rm d}\sg=\int_{{\bS}}\wt{B}(v-v_*,\og)\psi(v',v_*')\big|_{\og-{\rm rep.}}{\rm d}\og.\eeas
Thus any integral property which is valid for the $\og$-representation (or the $\sg$-representation)
is also valid for the $\sg$-representation (resp. the $\og$-representation).

In this paper we always use the $\sg$-representation (\ref{1.1}) and $B(z,\sg)$. Observe that the integrand (\ref{1.2}) is even in $\sg\in{\bS}$. This implies that
the collision kernel $B(z, \sg)$ in Eq.(BFD) can be replaced by its symmetrization
$\fr{1}{2}(B(z, \sg)+
B(z, -\sg)).$ In this paper we always assume that
$B(z,\sg)\equiv B(z,-\sg)$.

An important example is the Coulomb potential $\phi(|x|)={\rm  const.}|x|^{-1}$;
the corresponding collision kernel $B$ is given by the Rutherford formula
$B(z,\sg) = {\rm const.}\fr{1}{|z|^{3}\sin^{4}(\theta/2)}$; its symmetrization is
\be B(z,\sg) = {\rm const.}\fr{1+\cos^2(\theta)}{|z|^{3}\sin^{4}(\theta)}.
\lb{Coulomb}\ee
As pointed out by Bohm (pp.579, 555 in \cite{Bohm}) that the Rutherford formula ``has the unique property
that the exact classical theory, the exact quantum theory, and the Born approximation in the
quantum theory all yield the same scattering cross sections"(see also
(\ref{5.4*})). So far however it seems impossible to establish a weak form of Eq.(BFD) that can cover (\ref{Coulomb}) without angular cutoff. In this paper we assume a weaker angular cutoff for (\ref{Coulomb}) (see (\ref{Coul}) below).
On the other hand, it has been derived in\cite{weak-coupling},\cite{ESY} that in the weak-coupling regime the
collision kernel takes the following form (with $\cos\theta=\la z,\sg\ra/|z|, \theta\in[0,\pi]$)
\be B(z,\sg)=|z|\big(\wh{\phi}\big(|z\sin(\theta/2)|\big)-\wh{\phi}\big(|z\cos(\theta/2)|\big)
\big)^2\lb{weak-coupling}\ee
which is exactly computable in terms of (generalized) Fourier transform of the particle interaction potential
$ \phi(|x|)$:
$$\widehat{\phi}(|\xi|)=
\int_{{\bR}}\phi(|x|) e^{-{\rm i}\la\xi, x\ra}{\rm d}x,\quad \xi\in{\bR}. $$

Observe that the collision kernel $B(z,\sg)$ defined by (\ref{weak-coupling})
is even in $\sg\in {\bS}$.
For general case, since the integrand (\ref{1.2}) is even in $\sg\in{\bS}$,
the collision kernel $B(z, \sg)$ in Eq.(BFD) can be replaced by its symmetrization
$\fr{1}{2}(B(z, \sg)+
B(z, -\sg)).$ Thus we may always assume that
$B(z,\sg)\equiv B(z,-\sg)$.

Solutions of Eq.(BFD) are usually assumed having finite mass, momentum and kinetic energy, i.e.
$f(t,\cdot)\in L^1_2({\bR})$ for all $t\in [0,\infty)$, where $L^1_0({\bR})=L^1({\bR})$ and
$$L^1_s({\mR}^3)=\big\{f\in L^1({\mR}^3)\,\,|\,\,
\|f\|_{L^1_s}:=\int_{{\mR}^3}|f(v)|\la v\ra^{s}{\rm d}v<\infty\big\},\quad \la v\ra=(1+|v|^2)^{1/2},\quad s\in{\mR}.$$
In the following we denote
$$ \|f(t)\|=\|f(t,\cdot)\|,\quad  \int_{{\bR}}f(t){\rm d}v=\int_{{\bR}}f(t,v){\rm d}v,\quad \int_{{\bR}}g{\rm d}v=\int_{{\bR}}g(v){\rm d}v,\quad {\rm etc.}$$
Let $C^{m}_{b}([0,\infty)\times{\bR})$ be the space of functions $\vp\in C^m([0,\infty)\times{\bR})$
satisfying $\vp$ and all their partial derivatives of order $\le m$ are bounded on $[0,\infty)\times{\bR}$.

{\bf Definition of mild and weak solutions of Eq.(BFD).}
{\it Let $0\le f_0\in L^1_2({\bR})$ satisfy $0\le f_0\le 1$ on ${\bR}$. We say that a Lebesgue measurable function  $f: [0,\infty)\times{\bR} \to [0,1]$
satisfying
$\sup\limits_{t\in [0,\infty)}\|f(t)\|_{L^1_2}<\infty$ is a mild or weak solution of Eq.(BFD) with
initial datum $f(0,\cdot)=f_0$, if $f$ satisfies the following two conditions respectively:

Mild solution:\, There is a set $Z\subset {\bR}$ with $mes(Z)=0$ such that
for all $v\in {\bR}\setminus Z$,
\bes&& \int_{0}^{T}{\rm d}t
\int_{\bRS}B(v-v_*,\sg)\Pi_{{\rm F}}^{(\pm)}(f)(t,v,v_*,\sg){\rm d}\sg {\rm d}v_*
<\infty\qquad\forall\,
0<T<\infty,\lb{mild1}\\
&&
f(t,v)=f_0(v)+\int_{0}^{t}Q(f)(\tau, v){\rm d}\tau\qquad \forall\, t\in[0,\infty)
\dnumber\lb{mild2}\ees
where
\be\Pi_{{\rm F}}^{(+)}(f)=f'f'_*(1-f)(1-f_*),\quad
\Pi_{{\rm F}}^{(-)}(f)=ff_*(1-f')(1-f_*').\lb{pi1}\ee

Weak solution:\, For any $\vp\in C^{2}_{b}([0,\infty)\times{\bR})$ and any $0<T<\infty\,,$
\be\int_{0}^{T}{\rm d}t\int_{{\bRRS}}B(v-v_*,\sg)
|\Pi_{{\rm F}}(f)||\Dt\vp|\,{\rm d}\sg{\rm d}v_*{\rm d}v <\infty\,,\lb{weak1}\ee
$t\mapsto  \int_{{\bR}}f(t)\vp(t){\rm d}v$ is absolutely continuous on $[0,\infty)$, and
\be
\fr{{\rm d}}{{\rm d}t}\int_{{\bR}}f(t)\vp(t){\rm d}v
=\int_{{\bR}}f(t) {\p}_{t}\vp(t){\rm d}v- \fr{1}{4}
\int_{{\bR}}Q(f\,|\,\Dt\vp){\rm d}v\quad {\rm a.e.}\,\,\,t\in[0,\infty)\lb{weak2}\ee where
\bes&& Q(f\,|\,\Dt\vp)(t,v)
=\int_{{\bRS}}B(v-v_*,\sg) \Pi_{{\rm F}}(f)\Dt\vp{\rm d}\sg
{\rm d}v_*,\lb{weak3}\\
&&
\Dt\vp=\vp(t,v')+\vp(t,v_*')-\vp(t,v)-\vp(t,v_*)\,.\dnumber\lb{symdiff}\ees
If a (mild or weak) solution $f$  conserves the mass, momentum (or mean velocity) and kinetic energy, i.e. if
\be \int_{{\bR}} f(t,v)(1, v, |v|^2){\rm d}v=\int_{{\bR}} f_0(v)(1, v, |v|^2){\rm d}v\quad \forall\, t\in [0,\infty)\lb{conservation}\ee
then $f$ is called a conservative (mild or weak) solution of Eq.(BFD).
$\hfill\Box$}

In this paper we always assume that all initial data $f_0$ are non-vanished, i.e.
$\int_{{\bR}}f_0(v){\rm d}v>0$.

In order to prove the strong convergence to equilibrium, we use the following assumptions on the
collision kernel
$B(z,\sg)=B(|z|,\cos\theta)$:
there are Borel even functions $ b_*(\cdot), b^*(\cdot)$ in $(-1,1)$ and a Borel function
$\Phi_*$ in ${\mR}_{>0}$ satisfying
$$ \left\{\begin{array}
{ll}\displaystyle
0<b_*(t)\le b^*(t)\qquad \forall\, t\in (-1,1)\\
\\
\displaystyle
\Phi_*(r)>0\quad \forall\,r>0;\quad
\inf_{r\ge 1}\Phi_*(r)\ge 1,\quad \sup_{r>0}\Phi_*(r)<\infty
\end{array}\right.\eqno({\rm A1}) $$
such that
for all $(z,\sg)\in ({\bR}\setminus \{0\})\times {\bS}$ with $\cos\theta=\la z,\sg\ra/|z|$
the following (A2) or (A3) is satisfied:
$$ \left\{\begin{array}
{ll}\displaystyle |z|^{\gm}\Phi_*(|z|)b_*(\cos\theta)\le
B(z,\sg)\le(1+|z|^{\gm}) b^*(\cos\theta),\quad 0\le \gm\le 1
\\   \\ \displaystyle
 \int_{0}^{\pi}b^*(\cos\theta)
\sin\theta{\rm d}\theta<\infty\,.
\end{array}\right.\eqno({\rm A2}) $$
$$ \left\{\begin{array}
{ll}\displaystyle |z|^{\gm}\Phi_*(|z|)b_*(\cos\theta)\le
B(z,\sg)\le |z|^{\gm} b^*(\cos\theta),\quad -4\le \gm<0
\\   \\ \displaystyle
\int_{0}^{\pi}b^*(\cos\theta)
\sin^2(\theta){\rm d}\theta<\infty\,.
\end{array}\right.\quad\,\,\, \eqno({\rm A3})$$
Here the strict positivity in (A1) implies
$B(\cdot,\cdot)>0$ a.e. in ${\bRS}$ which
ensures that
every equilibrium of Eq.(BFD) must have the form (\ref{equilib1}) or (\ref{equilib2}).

In order to prove the strong convergence to equilibrium  we also
 assume that the lower bound functions $b_*(\cos\theta)$ in (A1)-(A3) are  {\it completely  positive} in $\cos^2(\theta)$, i.e.
$$
b_*(\cos\theta)=\sum_{n=0}^{\infty}a_n \cos^{2n}(\theta),\quad \theta\in(0,\pi),\quad a_n\ge 0,\,\,n=0,1,2,...\,.\eqno({\rm A4})$$
A weaker angular cutoff for the Rutherford formula (\ref{Coulomb}) and (\ref{5.4*}) is as follows
\be B(z,\sg)={\rm const}.|z|^{-3}\fr{1+\cos^2(\theta)}{\sin^{2p}(\theta)}\,\,\,
{\rm and}\,\,\, B(z,\sg)={\rm const}.|z|^{-3}\fr{\cos^2(\theta)}{\sin^{2p}(\theta)}, \quad 1<p<3/2. \lb{Coul}\ee
They satisfy the assumption (A3) for $\gm=-3$ with
$b_*(\cos\theta)=b^*(\cos\theta)={\rm const}.(1+\cos^2(\theta))
(1-\cos^2(\theta))^{-p}$ and $b_*(\cos\theta)=b^*(\cos\theta)={\rm const}.\cos^2(\theta)
(1-\cos^2(\theta))^{-p}$, which are completely positive in $\cos^2(\theta)$.
In Appendix we show some other typical collision kernels $B(z,\sg)$ of the form (\ref{weak-coupling}) that satisfy the assumptions \{(A1),(A2),(A4)\} and \{(A1),(A3),(A4)\} respectively.

It is easily seen that if $B(z,\sg)$ satisfies (A1),(A2) or $\{{\rm (A1)}, {\rm (A2)}_2\, {\rm and}\, {\rm (A3)}_1
 \,{\rm with} -3<\gm\le 1\}$\footnote{Here ${\rm (A2)}_2$ means the second condition in (A2).}, then the integrability condition
(\ref{mild1}) in the above definition is automatically satisfied, which is due to the $L^{\infty}$-bound $0\le f\le 1$ on ${\bR}\times[0,\infty)$ and $\sup\limits_{t\in [0,\infty)}\|f(t)\|_{L^1_2}<\infty$. And in that case, the existence and uniqueness of conservative mild solutions $f$ of Eq.(BFD) have been proven
and $f$ satisfy the moment estimates (\ref{moment4.1}) (for $0\le \gm\le 1$)
and the entropy identity (see e.g. \cite{Lu2001}, \cite{LW2003}):
\be S(f(t))-S(f_0)=\int_{0}^{t}D(f)(\tau){\rm d}\tau,\quad t\in [0,\infty)\lb{entropy-identity}\ee	
with \bes&&
S(f)=\int_{\bR}\big(-(1-f)\log(1-f)-f\log(f)\big){\rm d}v,\quad 0\le f\in L^1_2({\bR})\quad \lb{entropy}\\
&&
D(f)=\frac{1}{4}\int_{\bRRS}B(v-v_*,\sg)
\Gm\big(\Pi_{{\rm F}}^{(+)}(f),\Pi_{{\rm F}}^{(-)}(f)\big){\rm d}\sg {\rm d}v_*{\rm d}v\qquad
\dnumber\lb{entropy-D}\ees
where $\Gm(a,b)\in [0,\infty]$ is defined on ${\mR}_{\ge 0}^2$ by $\Gm(0,0)=0$ and
\be\Gm(a,b)=(a-b)\log\big(\fr{a}{b}\big)\,\,\,{\rm if}\,\,a,b>0;\,\,\,\Gm(a,b)=\infty\,\,\,{\rm if}\,\,\,
a>0=b\,\,{\rm or}\,\, a=0<b.\lb{Gamma}\ee
We note that the main results (existence, uniqueness, moment estimates, and the entropy identity) in \cite{Lu2001},\cite{LW2003} for Eq.(BFD) are actually proved for the $\og$-representation (\ref{1.8}) and for the
collision kernel $\wt{B}(z,\og)=b(\theta)|z|^{\gm}$\,($\theta=\arccos(|\la z,\og\ra|/|z|))$
of separated variable form. As mentioned above, these results hold also for the $\sg$-representation (\ref{1.4}); and the corresponding collision kernel $B$ can be relaxed to just satisfy the two-sided inequalities
in (A2)$_1$, (A3)$_1$ for $-3<\gm\le 1$ with the angular cutoff
(A2)$_2$.

The entropy identity (\ref{entropy-identity}) or entropy inequality (i.e. the left hand side $\ge $ the right hand side for (\ref{entropy-identity})) is
hardly proved to hold for weak solutions, because the function $(x_1,x_2,x_3,x_4)\mapsto \Gm\big(x_1x_2(1-x_3)(1-x_4), x_3x_4(1-x_1)(1-x_2)\big)$ is not convex
on $[0,1]^4$, and  it is very difficult to prove that
the $L^1([0,\infty)\times{\bRRS}, {\rm d}\nu)$-weak limiting function $q=q(t,v,v_*,\sg)$ of
a subsequence of $\{f^n f^n_*{f^n}'{f^n_*}'\}_{n=1}^{\infty}$
is equal to $ff_*f'f_*'$, where ${\rm d}\nu=e^{-t}e^{-|v-v_*|}{\rm d}t{\rm d}v{\rm d}v_*{\rm d}\sg$ and
$f^n$ are approximate solutions of Eq.(BFD) satisfying $f^n\rightharpoonup f$ weakly in $L^1([0,\infty)\times{\bR}, \nu)$. See \cite{Lu2008}. So far we can only prove the following
entropy inequality for some weak solutions $f$:
\be S(f(t))-S(f_0)\ge \int_{0}^{t}D(f,q)(\tau){\rm d}\tau,\quad t\in [0,\infty)\lb{weak-entropy1}\ee
where
\be D(f,q)(t)=\frac{1}{4}\int_{\bRRS}B(v-v_*,\sg)
\Gm\big(\Pi_{{\rm F}}^{(+)}(f,q),\Pi_{{\rm F}}^{(-)}(f,q)\big){\rm d}\sg {\rm d}v_*{\rm d}v
\lb{weak-entropy2}\ee
and the function $q=q(t,v,v_*,\sg)$ (mentioned above) is measurable on $[0,\infty)\times {\bRRS}$ satisfying
$0\le q\le 1$ on  $[0,\infty)\times {\bRRS}$
and for all $(t,v,v_*,\sg)\in [0,\infty)\times {\bRRS}$
\bes&& 0\le \Pi_{{\rm F}}^{(+)}(f,q)(t,v,v_*,\sg):=
f'f'_*(1-f-f_*)+q\le f'f'_*,\lb{weak-entropy3}\\
&& 0\le \Pi_{{\rm F}}^{(-)}(f,q)(t,v,v_*,\sg):=ff_*(1-f'-f'_*)+q\le ff_*. \dnumber
\lb{weak-entropy4}\ees
Although it is difficult to
prove that $q=ff_*f'f_*'$, the entropy inequality (\ref{weak-entropy1}) is enough for us to prove the strong convergence to equilibrium.

By entropy identity (\ref{entropy-identity})-(\ref{entropy-D})
and $B(\cdot,\cdot)>0$ a.e., an equilibrium $f(t, v)\equiv F(v)$ of Eq.(BFD) is
solution  of the functional equation $\Pi_{{\rm F}}(F)=0$ a.e. on ${\bRRS}$, i.e.
\be F'F_*'(1-F)(1- F_*)=FF_*(1-F')
(1-F_*') \ \ \ {\rm a.e.\,\,\,\, on}\quad {\bRRS}\lb{1.4}\ee
together with  the  physical  conditions
\be F\in L^1_2({\bR})\,\,\,\, {\rm and}\,\,\,\,
0\leq F\leq 1\,\,\,\,{\rm on}\,\,\,\,{\bR}.\lb{1.4*}\ee
It is proved in \cite{Lu2001} that for any $f_0\in L^1_2({\bR})$ satisfying $0\le f_0\le 1$ on ${\bR}$, it holds the inequality
\be \fr{M_2}{(M_0)^{5/3}}\ge \fr{3}{5}\Big(\fr{3}{4\pi}\Big)^{2/3}\lb{macg}\ee
where
\be M_0=\int_{{\bR}}f_0(v){\rm d}v,\quad M_2=\int_{{\bR}}|v-v_0|^2f_0(v){\rm d}v,\quad
v_0=\fr{1}{M_0}\int_{{\bR}}vf_0(v){\rm d}v,\lb{mac}\ee
and the corresponding equilibrium $F$ of Eq.(BFD) having the same mass, mean velocity and kinetic energy as $f_0$ is determined as follows:
\bes&&\fr{M_2}{(M_0)^{5/3}}> \fr{3}{5}\Big(\fr{3}{4\pi}\Big)^{2/3}
\quad \Longrightarrow\quad
F(v)=\fr{ae^{-b|v-v_0|^2}}{1+ae^{-b|v-v_0|^2}},\quad a>0, b>0;\lb{equilib1}\\
&&
\fr{M_2}{(M_0)^{5/3}}=\fr{3}{5}\Big(\fr{3}{4\pi}\Big)^{2/3}
\quad \Longrightarrow\quad
f_0(v)=F(v)={\bf 1}_{[0, R]}(|v-v_0|),\quad 0<R<\infty.\dnumber\lb{equilib2}\ees
According to the derivation in \cite{LW2003}, a main obstacle in proving the strong convergence to equilibrium is the nonnegativity problem of the quartic integral
\be \int_{\bRRS}B(v-v_*,\sg) hh_*h'h_*'{\rm d}\sg{\rm d}v{\rm d}v_*\,\ge 0\,? \lb{positivity}\ee
 for all $h\in L^1_2({\bR}),$ and $h$ will be often taken as $h(t,v)= f(t,v)-F(v)$. If there is no further restriction on $B(z,\sg)$, the quartic integral in (\ref{positivity}) can be negative for some $h\in L^1_2({\bR})$ (see Appendix for a counterexample). This is the only reason why the proof of strong convergence to equilibrium
given in \cite{LW2003} needs an additional assumption $a\le 1$ (which
is equivalent to a high temperature condition (\cite{LW2003})) which enables to avoid dealing
with the nonnegativity problem (\ref{positivity}), where $a$ is a coefficient of the equilibrium $F$
in (\ref{equilib1}). Our Proposition \ref{prop2.7} below
provides a class of collision kernels $B_0(z,\sg)$
(which are complete positive in $\cos^2(\theta)$ as shown in Appendix) that solves the nonnegativity problem (\ref{positivity}). This is the only reason
of using the assumption (A4). And the entropy
method also allows to use the lower bound kernel $|z|^{\gm}\Phi_*(|z|)b_*(\cos\theta)$ of $B(z,\sg)$, not necessarily requires $B(z,\sg)$ itself has the same property as
$|z|^{\gm}\Phi_*(|z|)b_*(\cos\theta)$.

Our main results of the paper are as follows:

\begin{theorem}\label{theorem1.1}
Suppose the collision kernel $B(z,\sg)$
satisfies (A1),(A2),(A4). Let $f_0\in L^1_2({\bR})$ satisfy
$0\le f_0\le 1$ on ${\bR}$,
and for the case $\gm=0$, assume in addition that $f_0\in L^1_s({\bR})$ for some $s>2$.
Let $f$ be the unique conservative mild solution of Eq.(BFD) with
the initial datum $f_0$ and let $F$ be the equilibrium of Eq.(BFD) having the same mass, momentum and kinetic energy as $f_0$. Then
\be \lim_{t\to\infty}\|f(t)-F\|_{L^1_2}=0.\lb{strong}\ee
\end{theorem}

\begin{theorem}\label{theorem1.2}
Suppose the collision kernel $B(z,\sg)$ satisfies (A1),(A3),(A4).
Let $f_0\in L^1_s({\bR})$ satisfy
$0\le f_0\le 1$ on ${\bR}$ with $s>\max\{2+|\gm|, 4\}$. Then there exists a conservative weak solution $f$ of Eq.(BFD) with the initial datum $f_0$ such that

(1) $f$ satisfies the entropy inequalities (\ref{weak-entropy1})-(\ref{weak-entropy4}) and the moment estimates:
\be\label{a-bounds}
	\|f(t)\|_{L^1_s}\leq C_s(1+t),\qquad \frac{1}{t}\int_0^t\|f(\tau)\|_{L^1_{s-|\gm|}}{\rm d}\tau\leq C_s,\qquad \forall\,t>0\ee
where the constant $0<C_s<\infty$ depends only on $b_*(\cdot), b^*(\cdot), \gamma, s, \|f_0\|_{L^1},\|f_0\|_{L^1_2}, $ and $\|f_0\|_{L^1_s}$.

(2) Let $F$ be the equilibrium of Eq.(BFD) having the same mass, momentum and kinetic energy as $f_0$.
Then
\be \lim_{t\to\infty}\fr{1}{t}\int_{0}^{t}\|f(\tau)-F\|_{L^1_2}{\rm d}\tau=0.\lb{averaging strong}\ee
\end{theorem}

\begin{remark}\lb{remark1} As mentioned above, Theorem \ref{theorem1.1} is
an improvement of \cite{LW2003} (with an addition assumption (A4)): the strong convergence (\ref{strong}) holds
for all temperature.  This theorem may be also established for weak
solutions $f$ under the assumptions (A1),(A2)$_1$,(A3)$_2$ and (A4);
but the more interesting thing is how to prove the
 convergence rate.
\end{remark}

\begin{remark}\lb{remark2} To our knowledge, Theorem \ref{theorem1.2} is the first result
concerning strong convergence to equilibrium for general initial data
of Eq.(BFD) for soft potentials ($-4\le \gm<0$). For the classical Boltzmann equation,
the averaging moment estimate in (\ref{a-bounds}) was first obtained in
\cite{D}  for $-1<\gm<0$ and
was generalized by \cite{CCL2009} to $-4\le \gm<0$, and in \cite{CCL2009} the time-averaged strong
 convergence to Maxwellian was proved for $-4\le \gm<0$ with and without
angular cutoff.
Here it is different from the classical case that the cubic integral
$\int_{{\bRS}}Bf_*f'f_*'{\rm d}\sg{\rm d}v_*$ in Eq.(BFD) does not have a useful property of separation of variables, and this is why we could only use a weaker angular cutoff (A3)$_2$
rather than the weakest angular cutoff $\int_{0}^{\pi}b^*(\cos\theta)\sin^3(\theta){\rm d}\theta<\infty$
as done in \cite{CCL2009}.
The averaged strong convergence (\ref{averaging strong}) in Theorem \ref{theorem1.2} may be improved as the usual strong convergence (\ref{strong});  the main difficulty in this regard perhaps lies in how to
prove $\liminf\limits_{T\to\infty}\inf\limits_{(t,v)\in[T,\infty)\times{\bR}}(1-f(t,v))>0$ for the case where the initial datum $f_0$ is not the second equilibrium in (\ref{equilib2}).
\end{remark}

The rest of the paper is organized as follows: In Section 2, under the assumptions (A1),(A3) we prove the moment estimates (\ref{a-bounds}) for every conservative weak solutions which satisfy the entropy inequalities (\ref{weak-entropy1})-(\ref{weak-entropy4}). The existence of such weak solutions is proved in Section 3. In Section 4 we use the entropy method, the nonnegativity shown in Proposition \ref{prop4.1}, and the regularity of the collision gain operator to prove the convergence results in Theorem
\ref{theorem1.1} and Theorem \ref{theorem1.2}.

\begin{center}\section{Moment Estimates}\end{center}

The moment estimates for mild solutions
under the assumptions (A1),(A2) has been proven in \cite{Lu2001}.
Here we prove moment estimates for weak solutions
under the assumptions (A1),(A3).

The following identity is well-known and is often used for investigating collision integrals:
$$\int_{{\bRRS}}B(v-v_*,\sg)\Psi(v',v_*',v,v_*){\rm d}\sg{\rm d}v_*{\rm d}v=
\int_{{\bRRS}}B(v-v_*,\sg)\Psi(v,v_*, v',v_*'){\rm d}\sg{\rm d}v_*{\rm d}v$$
for all nonnegative or integrable functions $\Psi$.

In this section for convenience of derivation we use a short notation:
$$\|f\|_s :=\|f\|_{L^1_s},\quad f\in L^1_s({\bR}).$$
We begin with several lemmas. The first lemma below is Proposition 2.1 in \cite{Lu2008}:

\begin{lemma}(\cite{Lu2008})\label{lemmA3.1}
Let $b, \Phi$ and $f$ be nonnegative Borel functions on $[0,\pi], {\mR}_{\ge 0}$ and
	${\bR}$ respectively.	Then for all $v\in {\bR}$,
 	$$\int_{\bRS}b(\theta)\Phi(|v-v_*|)f(v'){\rm d}\sg {\rm d}v_*=2\pi\int_{\bR}
 f(v_*){\rm d}v_*\int_{0}^{\pi}b(\theta)\frac{\sin(\theta)}{\sin^{3}(\theta/2)}
 \Phi\big(\fr{|v-v_*|}{\sin(\theta/2)}\big){\rm d}\theta,$$
$$\int_{\bRS}b(\theta)\Phi(|v-v_*|)f(v_*'){\rm d}\sg {\rm d}v_*=2\pi\int_{\bR}
	f(v_*){\rm d}v_*\int_{0}^{\pi}b(\theta)\frac{\sin(\theta)}{\cos^3(\theta/2)}
	\Phi\big(\fr{|v-v_*|}{\cos(\theta/2)}\big)
	{\rm d}\theta$$
	where in the left hand sides
	$\theta=\arccos(\la{\bf n},\sg\ra), {\bf n}=(v-v_*)/|v-v_*|$.
\end{lemma}

For any $\vp\in C^2({\bR})$ we use the following norm for Hesse matrix $\p^2\vp$:
$$
|\p^2\vp(v)|=\Big(\sum_{i=1}^3\sum_{j=1}^3|\p^2_{v_iv_j}\vp(v)|^2\Big)^{1/2}.$$

\begin{lemma}\lb{lemma3.2} (1)
 Let $\vp\in C^2(\bR), \Dt\vp=\vp(v')+\vp(v'_*)-\vp(v)-\vp(v_*).$ Then for all $\sg\in S^2,v,v_*\in \bR$
\be\lb{vp}
|\Dt\vp|\le\fr{1}{2}\Big(\sup_{|u|\le\sqrt{|v|^2+|v_*|^2}}|\p^2\vp(u)|\Big)|v-v_*|^2\sin\theta
\ee
where $\theta=\arccos(\la {\bf n},\sg\ra)\in [0,\pi], {\bf n}=(v-v_*)/|v-v_*|$.

(2) Let $b^*(\cdot)\ge 0$ satisfy (A3), $-5<\gm<0$, and let
$f\in L^1_2({\bR})$ satisfy $0\le f\le 1$ on ${\bR}$. Then for all $\vp\in C_b^2({\bR})$
\be
\int_{\bRRS}|v-v_*|^{\gm}b^*(\cos\theta)ff_*|\Dt\vp|{\rm d}\sg{\rm d}v_*{\rm d}v\le C_{\gm} A^*\|\p^2\vp\|_{L^{\infty}}\|f\|_{L^1}\big(1+\|f\|_{L^1_2}\big) \lb{qq3}\ee
where $A^*=2\pi\int_{0}^{\pi}b^*(\cos\theta)\sin^2(\theta){\rm d}\theta$,   $C_{\gm}=\max\{4,\fr{4\pi}{5-|\gm|}\}$.
\end{lemma}

\begin{proof} Part (1) is a partial result of Lemma 2.1 in \cite{CCL2009}.
For part (2), the left hand side of (\ref{qq3}) is less than
$$\fr{A^*}{2}\|\p^2\vp\|_{L^{\infty}}\int_{{\bRR}}f(v)f(v_*)|v-v_*|^{2-|\gm|}{\rm d}v{\rm d}v_*.$$
Since $0\le f\le 1$ and $-5<\gm<0$, this yields (\ref{qq3}).
\end{proof}
\\

\begin{lemma}\label{lemmA3.3}  Let $s>2, \vp_s(v)=\la v\ra^s,
\vp_{s,n}(v)=\vp_s(v)\zeta(|v|^2/n)$, where
$\zeta\in C^{\infty}([0,\infty))$ is defined by $\zeta(r)=1$ for all $r\in [0,1]; \zeta(r)=0$ for all
$r\in [2,\infty)$ and $0\le \zeta(r)\le 1$ for all $r\in [0,\infty)$.
 Then $\vp_{s,n}\in C^2_b({\bR})$ and (denoting $(s-4)_{+}=\max\{s-4,0\}$)
\bes&& \Dt\vp_s
	\le \fr{1}{2}s(s-2)\big(\la v\ra^2+\la v_*\ra^2\big)^{\fr{(s-4)_{+}}{2}}
	|v||v_*||v-v_*|^2\sin\theta,\lb{2.3}\\
&& |\Dt\vp_s|, |\Dt \vp_{s,n}|\le C_{s}(\la v\ra^2+\la v_*\ra^2)^{\fr{s-2}{2}}|v-v_*|^2\sin\theta
\dnumber \lb{2.4}\ees
where $0<C_s\le \fr{3}{2}(s^2+s+6)\max\{
1, \sup\limits_{r\ge 0}(1+r)|\fr{{\rm d}}{{\rm d}r}\zeta(r)|, \sup\limits_{r\ge 0}(1+r)r|\fr{{\rm d}^2}{{\rm d}r^2}\zeta(r)|\}$, which is independent of $n$.
\end{lemma}

\begin{proof} Let $\psi(r)=(1+r)^{s/2}$. Then
$\vp_s(v)=\psi(|v|^2)$ and using
$|v_*|^2-|v_*'|^2=|v'|^2-|v|^2$  we have
\be\Dt\vp_s=\Dt\psi(|\cdot|^2)=(|v'|^2-|v|^2)(|v|^2-|v_*'|^2)
\int_{0}^{1}\int_{0}^{1}
	\fr{{\rm d}^2}{{\rm d}r^2}\psi (R_{r,\tau}){\rm d}t{\rm d}\tau \lb{sym}\ee
where
$R_{t,\tau}=|v_*'|^2+t(|v'|^2-|v|^2)+\tau(|v|^2-|v_*'|^2)\,(t,\tau\in[0,1])$
which lies between \\ $\min\{|v|^2, |v_*|^2,|v'|^2,|v_*'|^2\}$ and
$\max\{|v|^2, |v_*|^2,|v'|^2,|v_*'|^2\}$.
Since $\fr{{\rm d}^2}{{\rm d}r^2}\psi(r)=\fr{s}{2}(\fr{s}{2}-1)(1+r)^{\fr{s}{2}-2}>0$ for
all $r\ge 0$, to prove (\ref{2.3}) we need only to consider the case
$(|v'|^2-|v|^2)(|v|^2-|v_*'|^2)>0$.
If $|v'|^2-|v|^2>0 $ and $|v|^2-|v_*'|^2>0$, then
$|v_*|^2-|v_*'|^2=|v'|^2-|v|^2>0$ and so
$|v_*|>|v_*'|$ and $|v|>|v_*'|$. Since
$|v_*-v_*'|=|v'-v|=|v-v_*|\sin(\theta/2),
|v-v_*'|=|v_*-v'|=|v-v_*|\cos(\theta/2)$, it follows that
\be (|v'|^2-|v|^2)(|v|^2-|v_*'|^2)=(|v_*|^2-|v_*'|^2)(|v|^2-|v_*'|^2)
\le 2|v||v_*||v-v_*|^2\sin(\theta).\lb{VV}\ee
Similarly if $|v'|^2-|v|^2<0 $ and $|v|^2-|v_*'|^2<0$, then
$|v'|<|v|$ and $|v_*|^2-|v'|^2=|v_*'|^2-|v|^2>0$
and so (\ref{VV}) still holds true.
From (\ref{sym})-(\ref{VV}), the inequality (\ref{2.3}) is easily deduced. The inequality (\ref{2.4}) is also easily
deduced by simple calculation using (\ref{vp}).
\end{proof}
\\

Let $B(z,\sg)=B(|z|,\cos\theta)$ satisfy (A3).
Define for all $\vp\in C^2(\bR)$ and $v,v_*\in{\bR}$,
\be
L[\Dt\vp](v,v_*)=\int_0^{\pi}B(|v-v_*|,\cos\theta)\sin(\theta)\Big(\int_{{\mS}^1(\textbf{n})}\Dt\vp{\rm d}\widetilde\sg\Big){\rm d}\theta
\ee
where ${\bf n}=(v-v_*)/|v-v_*|$ and ${\bf n}=(1,0,0)$ for $v=v_*$;
${\mS}^1({\bf n})=\{\wt{\sg}\in {\bS}\,\,|\,\,
\la \wt{\sg},{\bf n}\ra=0\}$, and ${\rm d}\wt{\sg}$ is the
length-measure element on ${\mS}^1({\bf n})$.

\begin{lemma}\lb{lemmA3.4} Let $\vp_s(v)=\la v\ra^s, L[\Dt\vp_s](v,v_*)$ be defined above. Then for any $s>2$ we have
	\begin{itemize}
		\item [(1)]
		If $-2\le\gm<0,$ then for any $\vep>0$
$$L[\Dt\vp_s](v,v_*)\le-c_s\big(\la v\ra^{s-|\gm|}+\la v_*\ra^{s-|\gm|}\big)+\vep C_s\big(\la v\ra^{s-|\gm|}\la v_*\ra^2+\la v_*\ra^{s-|\gm|}\la v\ra^2\big)+C_{s,\vep}\la v\ra^2\la v_*\ra^2.\quad $$
	\end{itemize}
	\begin{itemize}
		\item [(2)] If $-4\le\gm<-2$, then for any $\ld\ge1$
		\beas
		\begin{split}
			&{\bf 1}_{\{|v-v_*|>\ld\}}L[\Dt\vp_s](v,v_*)\\
			&\le-c_s\big(\la v\ra^{s-|\gm|}+\la v_*\ra^{s-|\gm|}\big)+
C_s\ld^{2-|\gm|}\big(\la v\ra^{s-|\gm|}\la v_*\ra^2+\la v_*\ra^{s-|\gm|}\la v\ra^2\big)+C_{s,\ld}\la v\ra^2\la v_*\ra^2.\quad
		\end{split}
		\eeas
	\end{itemize}
Here $0<c_s, C_s<\infty$ depends only on $\int_{0}^{\pi}b_*(\cos\theta)\sin^3(\theta){\rm d}\theta, \int_{0}^{\pi}b^*(\cos\theta)\sin^3(\theta){\rm d}\theta, \gm,$ and $s$.
\end{lemma}

\begin{proof} This is essentially the result of Lemma 4.2 in \cite{CCL2009}; the only difference
is that the lower bound function of $|z|$ in (A3) is $|z|^{\gm}\Phi_*(|z|)$ rather than
$(1+|z|)^{\gm}$ as used in \cite{CCL2009}. But the former can be derived by the latter:
let
$$B^{(1)}(z,\sg)=B(z,\sg)+|z|^{\gm}{\bf 1}_{\{|z|< 1\}}b_*(\cos\theta).$$
By $\gm<0$ we have
$(1+|z|)^{\gm}\le |z|^{\gm}\Phi_*(|z|)+|z|^{\gm}{\bf 1}_{\{|z|< 1\}}$ and so
$$(1+|z|)^{\gm}b_*(\cos\theta)
\le B^{(1)}(z,\sg)\le |z|^{\gm}\big(b^*(\cos\theta)+b_*(\cos\theta)\big).$$
So $B^{(1)}(z,\sg)$ satisfies the assumption in \cite{CCL2009}. Let
$L[\Dt\vp_s], L^{(1)}[\Dt\vp_s]$ correspond to $B(z,\sg), B^{(1)}(z,\sg)$
respectively. Then according to \cite{CCL2009}, part (1) and part (2)
hold for $ L^{(1)}[\Dt\vp_s]$.

Now suppose $-2\le \gm<0$. Let $L^{(2)}[\Dt\vp_s]$ correspond to $B^{(2)}(z,\sg):=|z|^{\gm}{\bf 1}_{\{|z|<1\}}b_*(\cos\theta).$
Then $L[\Dt\vp_s](v,v_*)=L^{(1)}[\Dt\vp_s](v,v_*)-L^{(2)}[\Dt\vp_s](v,v_*)$
and using (\ref{2.4}) we have
\be
|L^{(2)}[\Dt\vp_s](v,v_*)|\le C_s(\la v\ra^2+\la v_*\ra^2)^{\fr{s-2}{2}}|v-v_*|^{2-|\gm|}{\bf 1}_{\{|v-v_*|<1\}}\lb{2.7}\ee
with a different constant $C_s$.
It is easily seen that $|v-v_*|<1\Longrightarrow
\la v\ra^2+\la v_*\ra^2< 4(\la v\ra\wedge \la v_*\ra)^2$. This
together with $|\gm|\le 2$ and $\la \cdot\ra\ge 1$ implies that the right hand side of
(\ref{2.7}) is less than $
C_s(\la v\ra\wedge \la v_*\ra)^{s-2}$. Applying the inequality
$a^{k-2}\le \vep a^k+(1+\vep^{-\fr{k-4}{2}}) a^2, a\ge 1, k\in {\mR},$ to
$a=\la v\ra\wedge \la v_*\ra$ and $k=s-|\gm|$ we deduce
(with a different $C_s$) that
$C_s(\la v\ra\wedge \la v_*\ra)^{s-2}=C_s(\la v\ra\wedge \la v_*\ra)^{s-|\gm|-2}
(\la v\ra\wedge \la v_*\ra)^{|\gm|}
\le C_s\vep \la v\ra ^{s-|\gm|}\la v_*\ra^2+C_{s,\vep}\la v\ra^2\la v_*\ra^2.$
Thus
$L[\Dt\vp_s](v,v_*)\le L^{(1)}[\Dt\vp_s](v,v_*)+C_s\vep \la v\ra ^{s-|\gm|}\la v_*\ra^2+C_{s,\vep}\la v\ra^2\la v_*\ra^2$
and so part (1) holds for $L[\Dt\vp_s]$.

Next suppose $-4\le \gm<-2$. Then for any $\ld \ge 1$ we have
${\bf 1}_{\{|v-v_*|>\ld\}}L[\Dt\vp_s](v,v_*)=
{\bf 1}_{\{|v-v_*|>\ld\}}L^{(1)}[\Dt\vp_s](v,v_*)$ and so part (2) holds also
for $L[\Dt\vp_s]$.
\end{proof}
\\

\begin{lemma}\label{lemmA3.5} Let
	$b^*(\cdot)$ be the positive Borel function in (A1),(A3), let $k\ge 0, \beta\in{\mR}$ and let
	$f\in L^1_s({\bR})$ satisfy $0\le f\le 1$ on ${\bR}$ with $s\ge \max\{k, \beta, k+\beta\}$.
Then the following inequalities hold for all $v\in {\bR}$:
	
(1) Suppose $\beta\ge 0$.  Then
	\be\int_{{\bRS}}b^*(\cos\theta)\sin(\theta)f'f_*'\min\{\la v'\ra^k, \la v_*'\ra^k\}
	|v-v_*|^{\beta}
	{\rm d}\sg{\rm d}v_*
\le C_{\beta}\big(\|f\|_k\la v\ra^{\beta}+\|f\|_{k+\beta}\big).\lb{cub1} \ee
	And for any $0<\vep<\pi/2$
\be
	\int_{{\bRS}}b^*(\cos\theta)\sin(\theta)f'f_*'
	|v-v_*|^{\beta}
	{\rm d}\sg{\rm d}v_*
\le \fr{C_{\beta}}{\sin^3(\vep/2)}\|f\|_{\beta}
+C_{\beta}A^*(\vep)\|f\|_0\la v\ra^\beta\lb{cub2}\ee
where
\be A^*(\vep)=2\pi\int_{0}^{\vep}b^*(\cos\theta)\sin^2(\theta){\rm d}\theta+2\pi\int_{\pi-\vep}^{\pi}b^*(\cos\theta)\sin^2(\theta){\rm d}\theta.\lb{AE}\ee
	
(2) Suppose $\beta<0$. Then for any $\ld\ge 1$,
\be
\int_{{\bRS}}b^*(\cos\theta)\sin(\theta){\bf 1}_{\{|v-v_*|>\ld\}}f'f_*'\min\{\la v'\ra^k, \la v_*'\ra^k\}|v-v_*|^{\beta}{\rm d}\sg{\rm d}v_*
\le C_{\beta} \ld ^{\beta}\|f\|_k,\lb{cub3}\ee
\be\int_{{\bRS}}b^*(\cos\theta)\sin(\theta){\bf 1}_{\{|v-v_*|>\ld\}}f'f_*'\min\{\la v'\ra^k, \la v_*'\ra^k\}
|v-v_*|^{\beta-1}{\rm d}\sg{\rm d}v_*
\le C_{\beta} \ld^{\beta}\|f\|_{k-1}\la v\ra.\lb{cub4}\ee
The constants $C_\beta$ depend only on $\beta$ and $A^*:=2\pi\int_{0}^{\pi}b^*(\cos\theta)\sin^2(\theta){\rm d}\theta$.	
\end{lemma}

\begin{proof} Let $\Phi(r)$ be one of the functions
$r^{\beta}, r^{\beta}{\bf 1}_{\{r\ge \ld\}}, r^{\beta-1}{\bf 1}_{\{r\ge \ld\}}$ for $r\in{\mR}_{>0}.$
By $0\le f\le 1$ we have $f'f_*'\min\{\la v'\ra^k, \la v_*'\ra^k\}
\le {\bf 1}_{\{0\le \theta< \pi/2\}}f_*'\la v_*'\ra^k+{\bf 1}_{\{\pi/2<\theta\le \pi\}}f'\la v'\ra^k
$ and so using Lemma \ref{lemmA3.1} we have
\beas I(f,\Phi,k):&=&
\int_{{\bRS}}b^*(\cos\theta)\sin(\theta)f'f_*'\min\{\la v'\ra^k, \la v_*'\ra^k\}
\Phi(|v-v_*|){\rm d}\sg{\rm d}v_*\\
&\le& 2\pi\int_{{\bR}}f(v_*)\la v_*\ra^k
	{\rm d}v_*\int_{0}^{\pi/2}b^*(\cos\theta)
	\fr{\sin^2(\theta)}{\cos^{3}(\theta/2)}
	\Phi\big(\fr{|v-v_*|}{\cos(\theta/2)}\big)
	{\rm d}\theta
	\\
&+&2\pi\int_{{\bR}}f(v_*)\la v_*\ra^k
	{\rm d}v_*\int_{\pi/2}^{\pi}b^*(\cos\theta)
	\fr{\sin^2(\theta)}{\sin^{3}(\theta/2)}
	\Phi\big(\fr{|v-v_*|}{\sin(\theta/2)}\big)
	{\rm d}\theta.\eeas
	
	(1): Suppose $\beta\ge 0$. Then taking $\Phi(r)=r^{\beta}$ we obtain
\beas I(f,\Phi,k) \le
	2^{(3+\beta)/2}A^*\int_{{\bR}}f(v_*)\la v_*\ra^k
	|v-v_*|^{\beta}{\rm d}v_*.\eeas
Since
$\la v_*\ra^k|v-v_*|^{\beta}\le
2^{\beta}
	\big(\la v\ra^{\beta}\la v_*\ra^k+ \la v_*\ra^{k+\beta}\big)$,
this gives (\ref{cub1}).
Next  take any $0<\vep<\pi/2$ and recall that
$\theta=\arccos(\la {\bf n},\sg\ra), {\bf n}=(v-v_*)/|v-v_*|.$
We have
\beas&&
	\int_{{\bRS}}b^*(\cos\theta)\sin(\theta)f'f_*'
	|v-v_*|^{\beta}
	{\rm d}\sg{\rm d}v_*\\
	&&=
	\int_{{\bRS}}b^*(\cos\theta)\sin(\theta){\bf 1}_{\{\vep\le \theta\le \pi-\vep\}}f'f_*'
	|v-v_*|^{\beta}
	{\rm d}\sg{\rm d}v_*\\
	&&+
	\int_{{\bRS}}b^*(\cos\theta)\sin(\theta)\big({\bf 1}_{\{0\le \theta<\vep\}}+
{\bf 1}_{\{\pi-\vep<\theta\le \pi\}}
\big)f'f_*'
	|v-v_*|^{\beta}
	{\rm d}\sg{\rm d}v_*\\
	&&=: I_1+I_2.\eeas
	For $I_1$, using $0\le f\le 1$ again we have
$f'f_*'|v-v_*|^{\beta}
	\le  2^{\beta} \big(f'\la v'\ra^{\beta}+f_*'\la v_*'\ra^{\beta}\big)$
and so
\beas I_1&\le& 2^{\beta}2\pi\int_{{\bR}}f(v_*)\la v_*\ra^{\beta}
{\rm d}v_*
\int_{\vep}^{\pi-\vep}b^*(\cos\theta)
\fr{\sin^2(\theta)}{\sin^{3}(\theta/2)}
{\rm d}\theta
\\
&+&2^{\beta}2\pi\int_{{\bR}}f(v_*)\la v_*\ra^{\beta}
{\rm d}v_*
\int_{\vep}^{\pi-\vep}b^*(\cos\theta)
\fr{\sin^2(\theta)}{\cos^{3}(\theta/2)}
{\rm d}\theta\\
&\le & 2^{\beta+1}A^*\fr{1}{\sin^3(\vep/2)}\|f\|_{\beta}
.\eeas
For $I_2$, again using $({\bf 1}_{\{0\le \theta<\vep\}}+
{\bf 1}_{\{\pi-\vep<\theta\le \pi\}}
)f'f_*'\le {\bf 1}_{\{0\le \theta<\vep\}}f_*'+{\bf 1}_{\{\pi-\vep<\theta\le \pi\}}f'
$ and Lemma \ref{lemmA3.1} we have
\beas
I_2 &\le&
2\pi\int_{{\bR}}f(v_*)|v-v_*|^{\beta}{\rm d}v_*
	\int_{0}^{\vep}b^*(\cos\theta)\fr{\sin^2(\theta)}{\cos^{3+\beta}(\theta/2)}
	{\rm d}\theta
	\\
&+& 2\pi\int_{{\bR}}f(v_*)|v-v_*|^{\beta}{\rm d}v_*
	\int_{\pi-\vep}^{\pi}b^*(\cos\theta)\fr{\sin^2(\theta)}{\sin^{3+\beta}(\theta/2)}
	{\rm d}\theta
	\\
&\le &2^{(3+\beta)/2}A^*(\vep)\int_{{\bR}}f(v_*)|v-v_*|^{\beta}{\rm d}v_*.
\eeas
Then using $|v-v_*|^{\beta}\le 2^{\beta} \big(\la v\ra^{\beta}+\la v_*\ra^{\beta}\big)$ gives (\ref{cub2}).

(2): Suppose $\beta<0$.  Taking $\Phi(r)= r^{\beta}{\bf 1}_{\{r\ge \ld\}}$ we obtain
\beas I(f,\Phi,k)
\le 2^{(3+\beta)/2}A^*\int_{{\bR}}f(v_*)\la v_*\ra^k
|v-v_*|^{\beta}{\bf 1}_{\{\sqrt{2}|v-v_*|>\ld\}}{\rm d}v_*
\le 2^{3/2}A^*\ld^{\beta}\|f\|_{k}.\eeas
This proves (\ref{cub3}).
Taking $\Phi(r)= r^{\beta-1}{\bf 1}_{\{r\ge \ld\}}$ we also obtain
\bes I(f,\Phi,k)
&\le& C_{\beta}\int_{{\bR}}f(v_*)\la v_*\ra^k
|v-v_*|^{\beta-1}{\bf 1}_{\{\sqrt{2}|v-v_*|>\ld\}}{\rm d}v_*\nonumber\\
&\le& C_{\beta}\ld^{\beta}\int_{{\bR}}f(v_*)\fr{\la v_*\ra^k}{|v-v_*|}
	{\bf 1}_{\{|v-v_*|>1/\sqrt{2}\}}{\rm d}v_*=: C_{\beta}\ld^{\beta} J \lb{cubJ}\ees
where we used $\ld\ge 1$. Further estimates:
\beas J
&\le&
\int_{|v_*|<\la v\ra}f(v_*)\fr{\la v_*\ra^k}{|v-v_*|}
	{\bf 1}_{\{|v-v_*|>1/\sqrt{2}\}}{\rm d}v_*
+
\int_{|v_*|\ge \la v\ra}f(v_*)\fr{\la v_*\ra^k}{|v-v_*|}
	{\bf 1}_{\{|v-v_*|>1/\sqrt{2}\}}{\rm d}v_*\\
&\le& 2\|f\|_{k-1}\la v\ra +
\int_{|v_*|\ge \la v\ra}f(v_*)\fr{\la v_*\ra^k}{|v-v_*|}
	{\bf 1}_{\{|v-v_*|>1/\sqrt{2}\}}{\rm d}v_*=:J_1+J_2
.\eeas
For $J_2$, we consider a decomposition
\beas J_2 &=&
\int_{|v_*|\ge \la v\ra, |v-v_*|\le \fr{1}{2}|v_*|}f(v_*)\fr{\la v_*\ra^k}{|v-v_*|}
	{\bf 1}_{\{|v-v_*|>1/\sqrt{2}\}}{\rm d}v_*\\
&+&
\int_{|v_*|\ge \la v\ra, |v-v_*|> \fr{1}{2}|v_*|}f(v_*)\fr{\la v_*\ra^k}{|v-v_*|}
	{\bf 1}_{\{|v-v_*|>1/\sqrt{2}\}}{\rm d}v_*\\
&=:& J_{21}+J_{22}.\eeas
Note that $|v-v_*|\le \fr{1}{2}|v_*|$ implies
$\la v_*\ra\le 2\la v\ra$.
So $J_{21}\le 2\sqrt{2}\|f\|_{k-1}\la v\ra.$
Finally, $
|v_*|\ge \la v\ra$ and $|v-v_*|> \fr{1}{2}|v_*|$ imply
$|v-v_*|> \fr{1}{4}\la v_*\ra$. So
$J_{22}
\le 4\|f\|_{k-1}\le 4\|f\|_{k-1}\la v\ra.$
Thus we obtain
$J\le \big(2\sqrt{2} +6\big)\|f\|_{k-1}\la v\ra$. Connecting (\ref{cubJ}),
this proves (\ref{cub4}).
\end{proof}
\\

The following lemma deals with the cubic integral which is the main feature of quantum Boltzmann equations.

\begin{lemma}\label{lemmA3.6} Let $B=B(v-v_*,\sg)$ satisfy (A1),(A3) and let
$A^*(\vep)$ be defined in (\ref{AE}).
Let $s\ge\max\{2+|\gm|, 4\}, \vp_s(v)=\la v\ra^s$, and let
$f\in L^1_s({\bR})$ satisfy $0\le f\le 1$ on ${\bR}$. Then:
	
	If $-2\le \gm<0$, then  for any $0<\vep<\pi/2$
\be\int_{{\bRRS}}Bff'f_*'\Dt\vp_s {\rm}d\sg{\rm d}v_*{\rm d}v
\le  C_{s}\big(\vep+A^*(\vep)\big)\|f\|_2\|f\|_{s-|\gm|}+C_{s,\vep}\|f\|_2^2.
\lb{cub5}\ee

If $-4\le \gm<-2$, then for any $0<\vep<\pi/2, \ld\ge 1$
\be\int_{{\bRRS}}B_{\ld}ff'f_*'\Dt\vp_s {\rm d}\sg{\rm d}v_*{\rm d}v
\le C_{s}\big(\vep+A^*(\vep)+\ld^{\beta}\big)\|f\|_2\|f\|_{s-|\gm|}+
C_{s,\vep}\|f\|_2^2
\lb{cub6}\ee
where $B_{\ld}(v-v_*,\sg)=B(v-v_*,\sg){\bf 1}_{\{|v-v_*|>\ld\}},
\beta=2+\gm<0$ if $-3\le \gm<-2$; $\beta=3+\gm<0$ if  $-4\le \gm<-3$.
The constants $0<C_s<\infty$ depend only on $\gm,
A^*=2\pi\int_{0}^{\pi}b^*(\cos\theta)\sin^2(\theta){\rm d}\theta$, and
$s$, and
$0<C_{s,\vep}<\infty$ depend only on such constants $C_s$ and $\vep$.

\end{lemma}

\begin{proof} First of all using H\"{o}lder inequality and Young inequality
we have
\bes&&\alpha,\beta\ge 1, \alpha+\beta=s-|\gm|\,\Longrightarrow\,\|f\|_{\alpha}\|f\|_{\beta}\le \|f\|_1\|f\|_{s-|\gm|-1};\lb{2.15}\\
&&\dt>0\,\Longrightarrow\,
 \|f\|_{s-|\gm|-1}\le \dt \|f\|_{s-|\gm|}+C_{s,\dt}\|f\|_2\dnumber \lb{2.16}\ees
where $C_{s,\dt}$ depends only on $s,\gm,\dt$.
Next it is easily seen that
$\la v\ra^2, \la v_*\ra^2\le 2|v-v_*|^2+2\min\{\la v\ra^2, \la v_*\ra^2, \la v'\ra^2,
\la v_*'\ra^2\}$ and so
\be \big(\la v\ra^2+\la v_*\ra^2)^{\fr{s-4}{2}}|v_*|
 \le 2^{\fr{3}{2}(s-3)}\big(|v-v_*|^{s-3}+
\big(\min\{\la v\ra, \la v_*\ra, \la v'\ra,
\la v_*'\ra\}\big)^{s-3}\big).\lb{2.17}\ee
Let
$\chi(r)$ be one of the functions $1$ and ${\bf 1}_{\{r\ge \ld\}}$ for $
r\in{\mR}_{>0}$. Using inequalities (\ref{2.3}),
(\ref{2.17}) we have
\beas&&
 \int_{{\bRRS}}B(v-v_*,\sg)\chi(|v-v_*|)ff'f_*'\Dt\vp_s {\rm d}\sg{\rm d}v_*{\rm d}v
\\
&&\le C_s
\int_{{\bR}}f(v)\la v\ra {\rm d}v \int_{{\bRS}}b^*(\cos\theta))\sin(\theta)\chi(|v-v_*|)f'f_*'|v-v_*|^{s-|\gm|-1}{\rm d}\sg
{\rm d}v_*\\
&&+C_s
\int_{{\bR}}f(v)\la v\ra{\rm d}v
\int_{{\bRS}}b^*(\cos\theta))\sin(\theta)\chi(|v-v_*|)f'f_*'
\big(\min\{\la v'\ra,
\la v_*'\ra\}\big)^{s-3}|v-v_*|^{2-|\gm|}{\rm d}\sg
{\rm d}v_*
\\
&&=:I+J.\eeas
For the first term $I$, we take $\chi(|v-v_*|)\equiv 1$
and use  Lemma \ref{lemmA3.5}
with $\beta=s-|\gm|-1>0$  and use (\ref{2.16})
with $\dt=\vep \sin^3(\vep/2)$ to obtain
\beas I&\le& C_{s}\fr{1}{\sin^3(\vep/2)}\|f\|_1\|f\|_{s-|\gm|-1}
+C_{s}A^*(\vep)\|f\|_0\|f\|_{s-|\gm|}
\\
&\le& C_s\|f\|_1 \big(\vep+A^*(\vep)\big)\|f\|_{s-|\gm|}
+C_{s,\vep}\|f\|_1\|f\|_2.
\eeas
To estimate $J$, we consider three cases.

\noindent{\bf Case1}. $-2\le \gm<0$. Taking $\chi(|v-v_*|)\equiv 1$ and sing Lemma \ref{lemmA3.5}
with $k=s-3, \beta=2-|\gm|\ge 0$ we have
\beas J &\le&
 C_s\big( \|f\|_{s-3}\|f\|_{3-|\gm|}+\|f\|_1\|f\|_{s-|\gm|-1}\big)
\\
&\le& 2C_s\|f\|_1\|f\|_{s-|\gm|-1}
\le C_s\|f\|_1\vep\|f\|_{s-|\gm|}+C_{s,\vep}\|f\|_1\|f\|_2.
\eeas
Since $\|f\|_1\le \|f\|_2$, this proves (\ref{cub5}).

\noindent{\bf Case2.} $-3\le \gm<-2$. Taking $\chi(|v-v_*|)={\bf 1}_{\{|v-v_*|\ge \ld\}}$ and
using Lemma \ref{lemmA3.5} with $k=s-3, \beta=2-|\gm|<0$ we have
\beas&&
J\le C_s
\int_{{\bR}}f(v)\la v\ra\big(\ld^{\beta}\|f\|_{s-3}\big)
{\rm d}v
=C_s\ld^{\beta}\|f\|_1\|f\|_{s-3}\le C_s\ld^{\beta}\|f\|_1\|f\|_{s-|\gm|}.
\eeas

\noindent{\bf Case3.} $-4\le \gm<-3$. Taking $\chi(|v-v_*|)={\bf 1}_{\{|v-v_*|\ge \ld\}}$ and
using Lemma \ref{lemmA3.5} with $k=s-3, \beta=3-|\gm|<0$ we have
\beas
J\le
C_s\ld^{\beta}\|f\|_2\|f\|_{s-4}\le C_s\ld^{\beta}\|f\|_2\|f\|_{s-|\gm|}.
\eeas
Since $\|f\|_1\le \|f\|_2$, this together with the estimate for $I$ completes the proof of (\ref{cub6}).
\end{proof}
\vskip2mm

Having made enough preparations, we now prove moment estimates for
soft potentials ($-4\le \gm<0$).

\begin{proposition}\lb{prop2.7}
Let $B(z,\sg)$ satisfy (A1),(A3) and let $f_0\in L^1_s({\bR})$ satisfy $0\le f_0\le 1$ on ${\bR}$
with $s\ge \max\{2+|\gm|, 4\}$. Let
$f$ be a conservative weak solution of Eq.(BFD) with $f|_{t=0}=f_0$ satisfying
the entropy inequalities (\ref{weak-entropy1})-(\ref{weak-entropy4}).
Then
\be
\|f(t)\|_{L^1_s}\leq C_s(1+t),\qquad \frac{1}{t}\int_0^t
\|f(\tau)\|_{L^1_{s-|\gm|}}{\rm d}\tau\leq C_s,\qquad \forall\, t>0
	\lb{moment1}\ee
Here the constants $0<C_s<\infty$ depends only on $b_*(\cdot), b^*(\cdot), \gamma,\|f_0\|_{L^1}, \|f_0\|_{L^1_2},\|f_0\|_{L^1_s}$ and $s$.
\end{proposition}

\begin{proof}
To prove (\ref{moment1}), we need only
to prove that there are constants $0<c_s, C_s<\infty$ that
depend only on $\gm, b_*(\cdot), b^*(\cdot),\|f_0\|_{0}, \|f_0\|_{2},\|f_0\|_{s}$
and $s$ (and $c_s, C_s$ may have different value in different lines of derivation), such that $f(t,\cdot)\in L^1_s({\bR})$ for all $t\in[0,\infty)$ and
\be\|f(t)\|_s+c_s\int_{0}^{t}\|f(\tau)\|_{s-|\gm|}
{\rm d}\tau\le C_{s}(1+t)
\quad \forall\,t\in[0,\infty).\lb{moment2}\ee
\noindent {\bf Step1.} We prove that
\be\lb{39bd}
\sup_{t\in[0,T]}\|f(t)\|_s<\infty\qquad \forall\, 0<T<\infty.
\ee
Let $\vp_{s,n}(v)=\la v\ra^s\zeta(|v|^2/n)$ be defined in Lemma \ref{lemmA3.3}.
Then using that lemma and recalling that $0\le f(t,v)\le 1$ on $[0,\infty)\times {\bR}$
we have
\bes&&\int_{\bR}f(t,v)\vp_{s,n}(v){\rm d}v=\int_{\bR}f_0(v)\vp_{s,n}(v){\rm d}v-\frac{1}{4}\int_0^t{\rm d}\tau\int_{{\bR}}Q(f|\Dt \vp_{s,n})(\tau,v){\rm d}v\nonumber\\
&&\le \|f_0\|_s +C_s
\int_{0}^t{\rm d}\tau\int_{{\bRR}}ff_*(\la v\ra^2+\la v_*\ra^2)^{\fr{s-2}{2}}|v-v_*|^{2-|\gm|}
{\rm d}v{\rm d}v_*.\lb{approx}\ees
If $-2\le \gm<0$,
then
\be (\la v\ra^2+\la v_*\ra^2)^{\fr{s-2}{2}}|v-v_*|^{2-|\gm|}\le
2^{s/2}\big(\la v\ra^{s-|\gm|}+\la v_*\ra^{s-|\gm|}+\la v\ra^{s-2}\la v_*\ra^{2-|\gm|}+
\la v_*\ra^{s-2}\la v\ra^{2-|\gm|}\big).\lb{vv}\ee
So by
$\|f(t)\|_{2-|\gm|}\le \|f(t)\|_2=\|f_0\|_2$ we obtain
\beas\int_{\bR}f(t,v)\vp_{s,n}(v){\rm d}v\le
\|f_0\|_{s} +C_s \int_{{\bR}}f(\tau,v)\la v\ra^{s-|\gm|}{\rm d}v,\quad t\ge 0.\eeas
If $-4\le \gm<-2$, then using $0\le f\le 1$ and $\|f(\tau)\|_{0}=\|f_0\|_{0}$ we have
$\int_{{\bR}}f(\tau,v_*)|v-v_*|^{2-|\gm|}{\rm d}v_*
\le \fr{4\pi}{5-|\gm|}+\|f_0\|_{0}$
and so
$$\int_{\bR}f(t,v)\vp_{s,n}(v){\rm d}v\le
\|f_0\|_{s} +C_s \int_{{\bR}}f(\tau,v)\la v\ra^{s-2}{\rm d}v,\quad t\ge 0.$$
Let $\beta=\min\{|\gm|,2\}$, then
letting $n\to\infty$ we conclude from Fatou's Lemma that
\be\int_{{\bR}}f(t,v)\la v\ra^{s}{\rm d}v\le
\|f_0\|_{s} +C_s \int_{0}^{t}{\rm d}\tau\int_{{\bR}}f(\tau,v)\la v\ra^{s-\beta}{\rm d}v\qquad \forall\, t\ge 0.\lb{2.23}\ee
Note that the nonnegativity $f\ge 0$ ensures that the inequality (\ref{2.23}) holds for all $s>2$
satisfying $f_0\in L^1_s({\bR})$. Now let
$s\ge \max\{2+|\gm|, 4\}$ be given in the proposition.
Let $m\in {\mN}$ satisfy $2+m\beta\le s< 2+(m+1)\beta$ and let
$s_j=s-(m+1-j)\beta, j=0,1,..., m+1$.
Then  from (\ref{2.23}) we have
$$\int_{{\bR}}f(t,v)\la v\ra^{s_{j}}{\rm d}v\le
\|f_0\|_{s_j} +C_{s_{j}} \int_{0}^{t}{\rm d}\tau\int_{{\bR}}f(\tau,v)\la v\ra^{s_{j-1}}{\rm d}v\qquad \forall\, t\ge 0$$
for all $j=1,2,..., m+1$. Since $s_0<2$ so that $\|f(\tau)\|_{s_0}
<\|f(\tau)\|_2\equiv\|f_0\|_2$ (by conservation of mass and energy), it follows from the induction on $j$ that $f(t,\cdot)\in L^1_{s_j}({\bR})$ and
$\sup\limits_{t\in [0,T]}\|f(t)\|_{s_j}<\infty$ for all $0<T<\infty$ and all
$j=1,2,..., m+1$. In particular (\ref{39bd}) holds true.

\noindent{\bf Step2.} Let $\vp_s(v)=\la v\ra^s, \vp_{s,n}(v)=\vp_s(v)\zeta(|v|^2/n)$ be defined in
Lemma \ref{lemmA3.3}. Using Lemma \ref{lemmA3.3} we have
\beas B(v-v_*,\sg)|\Pi_{{\rm F}}(f)||\Dt\vp_{s,n}|
\le C_sb^*(\cos\theta)\sin(\theta)(f'f_*'+ff_*)(\la v\ra^2+\la v_*\ra^2)^{\fr{s-2}{2}}|v-v_*|^{2-|\gm|}.
\eeas
By {\bf Step1} and dominated convergence theorem we have
by letting $n\to\infty$ to the first equality in (\ref{approx})
that
\be\|f(t)\|_s=\|f_0\|_s-\frac{1}{4}\int_0^t{\rm d}\tau\int_{{\bR}}Q(f|\Dt \vp_s)(\tau,v){\rm d}v
\quad \forall\,t\in[0,\infty).\lb{absl}\ee
and
$$\int_{0}^{t}{\rm d}\tau
\int_{{\bRR}}ff_*(\la v\ra^2+\la v_*\ra^2)^{\fr{s-2}{2}}|v-v_*|^{2-|\gm|}{\rm d}v{\rm d}v_*<\infty\quad
\forall\, t\in [0,\infty).$$
Also by the usual change of variables we have
\be-\frac{1}{4}\int_{{\bR}}Q(f|\Dt \vp_s){\rm d}v=\fr{1}{2}
\int_{{\bRR}}ff_*L[\Dt\vp_s]{\rm d}v{\rm d}v_*+\int_{{\bRRS}}Bf'f_*'f\Dt\vp_s
{\rm d}\sg{\rm d}v{\rm d}v_*.\quad \ee

\noindent{\bf Step3.} Suppose $-2\le \gm<0$.
In this case, using Lemma \ref{lemmA3.4}, Lemma \ref{lemmA3.6}, and $\|f(t)\|_2=\|f_0\|_2$
we have
\beas&&
-\frac{1}{4}\int_{{\bR}}Q(f|\Dt \vp_s)(t,v){\rm d}v\\
&&\le
-c_s\|f_0\|_0\|f(t)\|_{s-|\gm|}+ C_{s}\big(\vep+A^*(\vep)\big)\|f_0\|_2\|f(t)\|_{s-|\gm|}+C_{s,\vep}\|f_0\|_2^2.
\eeas
By continuity of $[0,\pi/2]\ni \vep\mapsto A^{*}(\vep)$ and $A^{*}(0)=0$,
there is an $0<\vep<\pi/2$
which depends only on $A^*(\cdot)$ and  $\fr{c_s\|f_0\|_0}{2C_s\|f_0\|_2}$ such that
$\vep+A^*(\vep)\le \fr{c_s\|f_0\|_0}{2C_s\|f_0\|_2}$ so that
$$-c_s\|f_0\|_0\|f(t)\|_{s-|\gm|}+ C_{s}\big(\vep+A^*(\vep)\big)\|f_0\|_2\|f(t)\|_{s-|\gm|}
\le -\fr{1}{2}c_s\|f_0\|_0\|f(t)\|_{s-|\gm|}$$
hence (by (\ref{absl}))
\beas&&\|f(t)\|_s\le \|f_0\|_s-
\fr{1}{2}c_s\|f_0\|_0\int_{0}^{t}\|f(\tau)\|_{s-|\gm|}
{\rm d}\tau+C_{s,\vep}\|f_0\|_2^2 t
\qquad \forall\,t\in[0,\infty).\eeas
This proves (\ref{moment2}).

\noindent{\bf Step4.} Suppose $-4\le \gm<-2$.
For any $\ld\ge1,$ we split $B=B^\ld+B_\ld$ with
\be  B_\ld(z,\sg)={\bf 1}_{\{|z|>\ld\}}B(z,\sg),\quad
B^\ld(z,\sg)={\bf 1}_{\{|z|\le\ld\}}B(z,\sg)\ee
and let $Q(\cdot), Q_{\ld}(\cdot),L_{\ld}[\cdot], Q^{\ld}(\cdot)$ be the operators corresponding to the kernels $B(z,\sg),B_{\ld}(z,\sg)$ and $B^{\ld}(z,\sg)$ respectively.
Then $Q(\cdot)=Q_{\ld}(\cdot)+Q^{\ld}(\cdot)$
and so
\be
\|f(t)\|_s=\|f_0\|_s
-\frac{1}{4}\int_{0}^{t}{\rm d}\tau\int_{{\bR}}Q_{\ld}(f|\Dt \vp_s){\rm d}v
-\frac{1}{4}\int_{0}^{t}{\rm d}\tau\int_{{\bR}}Q^{\ld}(f|\Dt \vp_s){\rm d}v.\lb{split}\ee
For the integral about $Q_{\ld}(\cdot)$ we have by Lemma \ref{lemmA3.4} and Lemma \ref{lemmA3.6} and the conservation of mass and energy  that for any $0<\vep<1$
\beas&&-\frac{1}{4}\int_{{\bR}}Q_{\ld}(f|\Dt \vp_s){\rm d}v=\fr{1}{2}
\int_{{\bRR}}ff_*L_{\ld}[\Dt\vp_s]{\rm d}v{\rm d}v_*+\int_{{\bRRS}}B_{\ld}f'f_*'f\Dt\vp_s
{\rm d}\sg{\rm d}v{\rm d}v_*\\
&&\le -2c_s\|f_0\|_0\|f(\tau)\|_{s-|\gm|}+
C_{s}\big(\vep+A^*(\vep)+\ld^{-\alpha}\big)\|f_0\|_2\|f(\tau)\|_{s-|\gm|}+
C_{s,\ld,\vep}\|f_0\|_2^2
\eeas
where $\alpha=\min\{|\beta|, |\gm|-2\}>0$.
So
\beas&&-\frac{1}{4}\int_{0}^{t}{\rm d}\tau\int_{{\bR}}Q_{\ld}(f|\Dt \vp_s){\rm d}v
\le -c_s\|f_0\|_0\int_{0}^{t}\|f(\tau)\|_{s-|\gm|}{\rm d}\tau\\
&&+
C_{s}\big(\vep+A^*(\vep)+\ld^{-\alpha}\big)\|f_0\|_2\int_{0}^{t}\|f(\tau)\|_{s-|\gm|}{\rm d}\tau+
C_{s,\ld,\vep}\|f_0\|_2^2 t.\eeas
To estimate the integral about $Q^{\ld}(\cdot)$, we will use the entropy inequality (\ref{weak-entropy1}).  Recall that  $\Gm(a,b)= (a-b)\log(a/b)$
defined in (\ref{Gamma}) and apply the elementary inequality
$|a-b|\leq \fr{1}{2}(\sqrt{a}+\sqrt{b})\sqrt{\Gm(a,b)}, a,b\geq 0,$
to $a=\Pi_{{\rm F}}^{(+)}(f,q), b=\Pi_{{\rm F}}^{(-)}(f,q)$, and
notice that $\Pi_{{\rm F}}(f)=\Pi_{{\rm F}}^{(+)}(f,q)-\Pi_{{\rm F}}^{(-)}(f,q)$
and $\sqrt{\Pi_{{\rm F}}^{(+)}(f,q)}+\sqrt{\Pi_{{\rm F}}^{(-)}(f,q)}\le \sqrt{f'f_*'}+\sqrt{ff_*}$.
We have
\be |\Pi_{{\rm F}}(f)|
\le \fr{1}{2}\big(\sqrt{f'f_*'}+\sqrt{ff_*}\big)
\sqrt{\Gm\big(\Pi_{{\rm F}}^{(+)}(f,q),\Pi_{{\rm F}}^{(-)}(f,q)\big)}
\lb{pi6}\ee
so that by Cauchy-Schwartz inequality we obtain
\beas&&\Big|-\frac{1}{4}\int_{{\bR}}Q^{\ld}(f|\Dt \vp_s)(t,v){\rm d}v\Big|
\le \frac{1}{4}\int_{{\bRRS}}B^{\ld}|\Pi_{{\rm F}}(f)||\Dt\vp_s|
{\rm d}\sg{\rm d}v_*{\rm d}v
\\
&&\le \frac{1}{4}\Big(\fr{1}{2}\int_{{\bRRS}}B^{\ld}
(f'f_*'+ff_*)|\Dt\vp_s|^2
{\rm d}\sg{\rm d}v_*{\rm d}v\Big)^{1/2}\sqrt{D(f,q)(t)}\eeas
and
\beas&&\fr{1}{2}\int_{{\bRRS}}B^{\ld}(f'f_*'+ff_*)
|\Dt\vp_s|^2
{\rm d}\sg{\rm d}v_*{\rm d}v=
\int_{{\bRRS}}B^{\ld}ff_*
|\Dt\vp_s|^2
{\rm d}\sg{\rm d}v_*{\rm d}v
\\
&&\le C_s\int_{{\bRRS}}{\bf 1}_{\{|v-v_*|<\ld\}}b^*(\cos\theta)\sin^2(\theta)
ff_*(\la v\ra^2+\la v_*\ra^2)^{s-2}|v-v_*|^{4-|\gm|}
{\rm d}\sg{\rm d}v_*{\rm d}v\\
&&\le C_sA^{*}\ld^{4-|\gm|}\int_{{\bRR}}{\bf 1}_{\{|v-v_*|<\ld\}}
ff_*(\la v\ra^2+\la v_*\ra^2)^{s-2}{\rm d}v_*{\rm d}v
.\eeas
By $\ld\ge 1$ we see that $|v-v_*|\le \ld$ impels
$\la v\ra,\la v_*\ra\le 3\ld\min\{\la v\ra, \la v_*\ra\}$
so that $\la v\ra^2+\la v_*\ra^2\le 6\ld \la v\ra\la v_*\ra$
and so
\beas&&C_sA^{*}\ld^{4-|\gm|}\int_{{\bRR}}{\bf 1}_{\{|v-v_*|<\ld\}}
ff_*(\la v\ra^2+\la v_*\ra^2)^{s-2}{\rm d}v_*{\rm d}v
\\
&&\le C_sA^{*}\ld^{s+2-|\gm|}\int_{{\bRR}}{\bf 1}_{\{|v-v_*|<\ld\}}
ff_*\la v\ra^{s-2}\la v_*\ra^{s-2}
{\rm d}v_*{\rm d}v
\le C_sA^{*}\ld^{s+2-|\gm|}\|f(t)\|_{s-2}^2.\eeas
Thus we obtain
\beas&&
\Big|-\fr{1}{4}\int_{0}^{t}{\rm d}\tau\int_{{\bR}}Q^{\ld}(f|\Dt \vp_s)(t,v){\rm d}v\Big|
\le C_s\ld^{(s+2-|\gm|)/2}\int_{0}^{t}\|f(\tau)\|_{s-2}
\sqrt{D(f,q)(\tau)}{\rm d}\tau.\eeas
Next by Cauchy-Schwartz inequality and $s-4\le s-|\gm|$ we have
 \beas
 \|f(\tau)\|_{s-2}\le\sqrt{\|f(\tau)\|_s}\sqrt{\|f(\tau)\|_{s-4}}\le\sqrt{\|f(\tau)\|_s}\sqrt{\|f(\tau)\|_{s-|\gm|}}
 \eeas
Thus for any $0<\vep<1$
\beas\ld^{(s+2-|\gm|)/2}\|f(\tau)\|_{s-2}
\sqrt{D(f,q)(\tau)}
\le \vep\|f(\tau)\|_{s-|\gm|}+\fr{1}{4\vep}\ld^{s+2-|\gm|}\|f(\tau)\|_sD(f,q)(\tau)\eeas
and so
\beas&&
\Big|-\fr{1}{4}\int_{0}^{t}{\rm d}\tau\int_{{\bR}}Q^{\ld}(f|\Dt \vp_s)(t,v){\rm d}v\Big|
\\
&&\le C_s \vep\int_{0}^{t}\|f(\tau)\|_{s-|\gm|}
{\rm d}\tau+\fr{C_s}{\vep}\ld^{s+2-|\gm|}\int_{0}^{t}\|f(\tau)\|_sD(f,q)(\tau)
{\rm d}\tau.\eeas
Combining the above estimates we deduce
\beas&&
\|f(t)\|_s
\le \max\{C_{s,\ld,\vep,}\|f_0\|_2^2, \|f_0\|_s\}(1+t)
-c_s\|f_0\|_0\int_{0}^{t}\|f(\tau)\|_{s-|\gm|}{\rm d}\tau\\
&&+
C_{s}\big(\vep+A^*(\vep)+\ld^{-\alpha}\big)\|f_0\|_2\int_{0}^{t}\|f(\tau)\|_{s-|\gm|}{\rm d}\tau
+C_{s,\vep,\ld}\int_{0}^{t}\|f(\tau)\|_sD(f,q)(\tau)
{\rm d}\tau.
\eeas
Choose $0<\vep<1, \ld>1$ which depend only
on $ A^*(\cdot), c_s, C_s,\|f_0\|_0, \|f_0\|_2$ such that
$$C_{s}\big(\vep+A^*(\vep)+\ld^{-\alpha}\big)\|f_0\|_2
 \le \fr{1}{2}c_s\|f_0\|_0.$$
Thus we obtain (for different constants $0<c_s, C_s<\infty$)
$$
\|f(t)\|_s+c_s\int_{0}^{t}\|f(\tau)\|_{s-|\gm|}{\rm d}\tau
\le C_s(1+t)+C_s\int_{0}^{t}\|f(\tau)\|_sD(f,q)(\tau)
{\rm d}\tau\quad \forall\, t\ge 0$$
and then applying Gronwall lemma we conclude
$$\|f(t)\|_s+c_s\int_{0}^{t}\|f(\tau)\|_{s-|\gm|}{\rm d}\tau
\le C_s(1+t)\exp\Big(C_s\int_{0}^{t}D(f,q)(\tau)
{\rm d}\tau\Big),\quad t\ge 0.$$
Finally using $0\le f=f(t,v)\le 1$ and $y|\log(y)|\le \sqrt{y}$ for $0\le y\le 1$
we have
$-(1-f)\log(1-f)-f\log(f)=
(1-f)\log(1+f/(1-f))+
f\log(1/f)
\le f+f\log(1/f){\bf 1}_{\{0<f\le e^{-|v|^2}\}}
+f\log(1/f){\bf 1}_{\{e^{-|v|^2}<f\le 1\}}
\le f+\sqrt{f}{\bf 1}_{\{0<f\le e^{-|v|^2}\}}+f\log(e^{|v|^2})
\le f+e^{-|v|^2/2}+|v|^2f$.
So
\be S(f(t))\le \|f(t)\|_2+\int_{{\bR}}e^{-|v|^2/2}{\rm d}v=
\|f_0\|_2+C_0. \lb{SLeq}\ee
By the entropy inequality (\ref{weak-entropy1}), this gives
\beas\sup_{t\ge 0}\int_{0}^{t}D(f,q)(\tau)
{\rm d}\tau\le \sup_{t\ge 0}\big(S(f(t))-S(f_0)\big)\le \|f_0\|_2+C_0\eeas
and proves (\ref{moment2}).
\end{proof}

\begin{center}\section{Existence of Solutions}\end{center}

As mentioned in the Introduction, under the assumptions (A1),(A2) for the collision kernel $B$, it has been proven in
\cite{Lu2001}, \cite{LW2003} that
for any initial datum $f_0\in L^1_2({\bR})$ satisfying
$0\le f_0\le 1$ on ${\bR}$, and for the case $\gm=0$ assuming in addition that $f_0\in L^1_s({\bR})$ for some $s>2$, there exists a unique conservative mild solution $f$ of Eq.(BFD) satisfying $f|_{t=0}=f_0$,
and $f$ also satisfy the moment estimates (\ref{moment4.1}) and entropy identity (\ref{entropy-identity}).
Thus here we need only to consider the case where $B$ satisfies (A1),(A3) and prove
the existence of conservative weak solutions $f$ of Eq.(BFD)  satisfying
the entropy inequalities (\ref{weak-entropy1})-(\ref{weak-entropy4}).
First of all we note that the weak form (\ref{weak1}) of Eq.(BFD) is well-defined. In fact from part (2) of Lemma \ref{lemma3.2} one sees that the integral over
${\bRRS}$ in the weak form (\ref{weak1}) of Eq.(BFD) is well-controlled
for all $f\in L^1_2({\bR})$ satisfying $0\le f\le 1$ and for all $\vp\in C_b^2({\bR})$:
$$\int_{{\bRRS}}B(v-v_*,\sg)
|\Pi_{{\rm F}}(f)||\Dt\vp|\,{\rm d}\sg{\rm d}v_*{\rm d}v
\le  2C_{\gm} A^*\|\p^2\vp\|_{L^{\infty}}\|f\|_{L^1}\big(1+\|f\|_{L^1_2}\big).$$
Now let
$$ B_n(z,\sg)
=\min\big\{B(z,\sg),\, n\big\}$$
and let $\wt{b}_*(\cos\theta)=\min\{b_*(\cos\theta), 1\},
\wt{\Phi}_*(|z|)=
\min\{\Phi_*(|z|),|z|^{-\gm}\}$. We have
for all $(z,\sg)\in ({\bR}\setminus\{0\})\times {\bS}$,
$$|z|^{\gm}\wt{\Phi}_*(|z|)\wt{b}_*(\cos\theta)\le B_n(z,\sg)\le
b^*(\cos\theta)|z|^{\gm}\quad \forall\, n\in{\mN} $$
and
$$B_n(z,\sg)\nearrow\,B(z,\sg)\quad (n\to\infty).$$
Let $f_0\in L^1_{s}({\bR})$ satisfy $0\le f_0\le 1$ and $s\ge \max\{2+|\gm|, 4\}$.
For every $n\in{\mN}$, applying \cite{Lu2001} to the bounded kernel $B_n$,
there exists a unique conservative mild (hence weak) solution $f^n$ of Eq.(BFD) on ${\bR}\times [0,\infty)$ with $f^n|_{t=0}=f_0$ and $f^n$ satisfies the entropy identity (\ref{entropy-identity}) with the entropy dissipation $D=D_n$ corresponding to $B_n$. Since all $B_n$ satisfies the same condition as (A3),
it follows from
Proposition \ref{prop2.7} that $f^n$ have the moment estimates (\ref{moment1}) with
some constants $0<c_s,C_s<\infty$ that do not depend on $n$.
Thus the stability result in Theorem 1 in \cite{Lu2008}
can be applied to conclude that there exist a subsequence $\{f^{n_k}\}_{k=1}^{\infty}$ and a weak solution $f$ of Eq.(BFD) with the initial datum $f_0$, such that for every $t\in [0,\infty)$,
$f^{n_k}(t,\cdot)\rightharpoonup f(t,\cdot)$ weakly in $L^1({\bR})$, and $f$
satisfies the entropy inequalities (\ref{weak-entropy1})-(\ref{weak-entropy4}). That $f$ conserves the mass, momentum, and energy
follows from this weak convergence, $s>2$, and $\sup\limits_{n\ge 1}\sup\limits_{t\in [0,T]}\|f^n(t)\|_{L^1_s}<\infty\,(\forall\, 0<T<\infty)$.
It should be noted that the very soft condition $-5<\gm\le -3$ in \cite{Lu2008} for inhomogeneous
solution $f(t,x,v)$ of Eq.(BFD) on $[0,\infty)\times {\mathbb T}^3\times {\bR}$
(here ${\mathbb T}^3$ is a periodic box) is only used to avoid dealing with the
boundedness problem of the square integrals
$$ \sup_{n\ge 1}\int_{{\mathbb T}^3}\big(\|f^n(t,x,\cdot)\|_{L^1}\big)^2
{\rm d}x < \infty\, ?  $$
Since our conservative approximate solutions $f^n=f^n(t,v)$ do not depend on $x$,
there is no such problem, and the proofs for the weak stability and the existence of weak solutions of Eq.(BFD) are much easier than those in \cite{Lu2008}.

\begin{center}\end{center}
\begin{center}\section{Convergence to Equilibrium}\end{center}
	
In this section we prove Theorem \ref{theorem1.1} and  Theorem \ref{theorem1.2}. To do this we first
prove some useful properties. We will use short notations
$\|f\|_{L^p}=\|f\|_{L^p({\bR})}, \|f\|_{L^{\infty}}=
\|f\|_{L^{\infty}({\bR})}$ for $f\in L^p({\bR})$ or $L^{\infty}({\bR})$; and
$\|\Psi\|_{L^p}=\|\Psi\|_{L^p({\bRR})}, \|\Psi\|_{L^{\infty}}=
\|\Psi\|_{L^{\infty}({\bRR})}$ for  $\Psi\in L^p({\bRR})$ or
$L^{\infty}({\bRR})$.

\vskip2mm

\begin{proposition}\lb{prop4.1}  Let $B_0(z,\sg)=B_0(|z|, \cos\theta), z\in{\bR}\setminus \{0\},
\cos\theta =\la z,\sg\ra/|z|,\sg\in {\bS},$ and suppose that $B_0(z,\sg)$ is completely
positive in $\cos^2(\theta)$ which means that there are Borel functions $a_n(r)$
 on ${\mR}_{>0}$ such that for all $z\in {\bR}\setminus\{0\}, \,\theta\in(0,\pi),$
\be B_0(|z|,\cos\theta)=\sum_{n=0}^{\infty}
a_n(|z|)\cos^{2n}(\theta),\quad a_n(|z|)\ge 0,\,\, n=0,1,2,...\,.\lb{positive0}\ee
Then for any Borel function $h:{\bR}\to {\mR}$ satisfying
$$\int_{{\bRRS}}B_0(v-v_*,\sg)|h(v)h(v_*)h(v')h(v_*')|{\rm d}\sg{\rm d}v{\rm d}v_*<\infty$$
it holds
\be \int_{{\bRRS}}B_0(v-v_*,\sg)h(v)h(v_*)h(v')h(v_*'){\rm d}\sg{\rm d}v{\rm d}v_*\ge 0.\lb{positive1}\ee
\end{proposition}

\begin{proof}  By change of variables we have
\beas \int_{\bRRS}B_0h'h'_*hh_*{\rm d}\sg {\rm d}v{\rm d}v_*
=\frac{1}{8}\int_0^\infty r^2{\rm d}r\int_{\mR^3} {\rm d}u\int_{\bSS}B_0(r,\la \og,\sg\ra)\wt{h}_{u,r}(\og)\wt{h}_{u,r}(\sg){\rm d}\og {\rm d}\sg\eeas
where
$$\wt{h}_{u,r}(\sg)=h\big(\frac{u+r\sg}{2})h(\frac{u-r\sg}{2}\big).$$
Thus, to prove (\ref{positive1}), we need only to prove that for any $r\in{\mR}_{>0}$, if a Borel function $h:{\bS}\to {\mR}$ satisfies
$$\int_{{\bSS}}B_0(r,\la \og,\sg\ra)|h(\og)h(\sg)|{\rm d}\og{\rm d}\sg<\infty$$
then
\be \int_{{\bSS}}B_0(r,\la \og,\sg\ra)h(\og)h(\sg){\rm d}\og{\rm d}\sg\ge 0.\lb{positive2}\ee
Let
$$B_m(r,t)=a_0(r)+\sum\limits_{n=1}^{m}a_n(r)t^{2n},\quad m\in{\mN}.$$
Then $|B_m(r,\la \og,\sg\ra)h(\og)h(\sg)|\le B_0(r,\la \og,\sg\ra)|h(\og)h(\sg)|$ and so by
dominated convergence theorem we have
$$\int_{\bSS}B_0(r,\la \og,\sg\ra)h(\og)h(\sg){\rm d}\og {\rm d}\sg=\lim_{m\to\infty}
\int_{\bSS}B_m(r,\la \og,\sg\ra)h(\og)h(\sg){\rm d}\og {\rm d}\sg.$$
Therefore to prove (\ref{positive2}) we need only to prove that for any $m\in{\mN}$,
\be\int_{{\bSS}}B_m(r,\la \og,\sg\ra)h(\og)h(\sg){\rm d}\og{\rm d}\sg\ge 0.\lb{positive3}\ee
To do this we first assume that all functions $a_n(r)$ are bounded on ${\mR}_{\ge 0}$ and $h$ is also
bounded on ${\bS}$ so that there is no problem of integrability in the following calculation.
Write  $\og=(\og_1, \og_2,\og_3), \sg=(\sg_1, \sg_2,\sg_3)$. Then
$\la \og,\sg\ra=\og_1\sg_1+\og_2\sg_2+\og_3\sg_3$,
\beas&&(\la \og,\sg\ra)^{2n}=
\sum_{i_1=1}^3\sum_{i_2=1}^{3}\cdots\sum_{i_{2n}=1}^3\og_{i_1}\og_{i_2}\cdots \og_{i_{2n}}
\sg_{i_1}\sg_{i_2}\cdots \sg_{i_{2n}},\\
&&
\int_{{\bSS}}B_m(r,\la \og,\sg\ra)h(\og)h(\sg){\rm d}\og{\rm d}\sg
=a_0(r)\int_{{\bSS}}h(\og)h(\sg){\rm d}\og{\rm d}\sg\\
&&
+\sum_{n=1}^m a_{n}(r)
\sum_{i_1=1}^3\sum_{i_2=1}^{3}\cdots\sum_{i_{2n}=1}^3\int_{{\bSS}}\og_{i_1}\og_{i_2}\cdots \og_{i_{2n}}
\sg_{i_1}\sg_{i_2}\cdots \sg_{i_{2n}}h(\og)h(\sg){\rm d}\og{\rm d}\sg
\\
&&=a_0(r)\Big(\int_{{\bS}}h(\sg){\rm d}\sg\Big)^2
+\sum_{n=1}^m a_{n}(r)
\sum_{i_1=1}^3\sum_{i_2=1}^3\cdots\sum_{i_{2n}=1}^3\Big(\int_{{\bS}}\sg_{i_1}\sg_{i_2}\cdots \sg_{i_{2n}}h(\sg){\rm d}\sg
\Big)^2
\ge 0.\eeas
Next for general case, for every $m\in{\mN}$, we use approximation
 $$B_m^k(r,\la \og,\sg\ra)=\fr{a_0(r)}{1+\fr{1}{k}a_0(r)}+
 \sum\limits_{n=1}^{m}\fr{a_n(r)}{1+\fr{1}{k}a_n(r)}\la \og,\sg\ra^{2n},\quad k\in{\mN}.$$
We have $0\le B_m^k(r,\la \og,\sg\ra)\le B_m(r,\la \og,\sg\ra)$ and
 $$\int_{{\bSS}}B_m^k(r,\la \og,\sg\ra)\fr{h(\og)}{1+\fr{1}{k}|h(\og)|}\cdot\fr{h(\sg)}{1+\fr{1}{k}|h(\sg)|}{\rm d}\og{\rm d}\sg\ge 0.$$
Letting $k\to\infty$ we obtain (\ref{positive3}) by dominated convergence theorem.
\end{proof}
\\

Let $B(z,\sg)$ be the collision kernel stated in Theorem \ref{theorem1.1} and
Theorem \ref{theorem1.2} respectively, and let
$b_*(\cos\theta)$ be the corresponding lower bound function in the assumptions (A1)-(A4).
In (A4), let us replace $a_n$ with $\underline{a}_n=\min\{a_n, 2^{-n}\}, n=0,1,2,...,$.
Then the new function is bounded and is still completely positive in $\cos^2(\theta)$:
$$\underline{b}_*(\cos\theta):=\sum_{n=0}^{\infty}\underline{a}_n\cos^{2n}(\theta)\le \min\{ b_*(\cos\theta), 2\},\quad \theta\in [0,\pi].$$
 In the following we use a smaller kernel
\be B_0(z,\sg)=\fr{\min\{|z|^{\gm}\Phi_*(|z|), 1\}}{1+|z|^6}\underline{b}_*(\cos\theta),\quad
-4\le \gm\le 1\lb{smaller}\ee
which is of course completely
positive in $\cos^2(\theta)$ defined in Proposition \ref{prop4.1}, and we have
$$B_0(z,\sg)\le B(z,\sg)\quad \forall\, (z,\sg)\in{\bRS}.$$
Let
\be
\underline{A}_*=2\pi\int_{0}^{\pi}\underline{b}_{*}(\cos\theta)\sin(\theta){\rm d}\theta,\quad
A_0(|z|)=\int_{{\bS}}B_0(z,\sg){\rm d}\sg=\underline{A}_*
\fr{\min\{|z|^{\gm}\Phi_*(|z|), 1\}}{1+|z|^6}.\dnumber \lb{A_0}\ee
The gain term in the collision operator
corresponds in a natural way to the linear operators ${\cal Q}^+, {\cal Q}^+(\cdot|F) : L^{\infty}\cap L^1(\mR^3\times\mR^3)\to L^{\infty}\cap L^1(\bR)$:
\bes&& {\cal Q}^+(\Psi)(v)=\int_{\bRS}B_0(v-v_*,\sg)\Psi(v',v'_*){\rm d}\sg {\rm d}v_*,\lb{calQ}\\
&&{\cal Q}^+(\Psi| F)(v)=\int_{\bRS}B_0(v-v_*,\sg)\Psi(v',v'_*)F(v_*){\rm d}\sg {\rm d}v_*
\dnumber\lb{calQF}\ees
where $F$ is any bounded Lebesgue measurable function and will be taken a smooth equilibrium (\ref{equilib1}) of Eq.(BFD). For convenience of derivation we assume  that
\be F\in L^{\infty}({\bR}),\quad \|F\|_{L^{\infty}}\le 1.\lb{FF}\ee
For any $f,g\in L^{\infty}\cap L^1({\bR})$ define $(f\otimes g)(v,v_*)=f(v)g(v_*)$. Since $\sg\mapsto B_0(z,\sg)$ is even, this implies that
$${\cal Q}^{+}(f\otimes g)={\cal Q}^{+}(g\otimes f),\quad
{\cal Q}^{+}(f\otimes g| F)={\cal Q}^{+}(g\otimes f| F).$$
It is easily deduced that for all $\Psi\in L^{\infty}\cap L^1(\mR^3\times\mR^3)$,
\bes&& \|{\cal Q}^+(\Psi)\|_{L^{\infty}}, \|{\cal Q}^+(\Psi|F)\|_{L^{\infty}}
\le \|A_0(|\cdot|)\|_{L^1}\|\Psi\|_{L^{\infty}},\lb{QF1}\\
&&\|{\cal Q}^+(\Psi)\|_{L^1}, \|{\cal Q}^+(\Psi|F)\|_{L^1}\le
\underline{A}_*\|\Psi\|_{L^1},\dnumber\lb{QF2}\\
&& \|{\cal Q}^+(\Psi)\|_{L^2}, \|{\cal Q}^+(\Psi|F)\|_{L^2}
\le \|A_0(|\cdot|)\|_{L^2}\|\Psi\|_{L^2}.\dnumber\lb{QF3}\ees
As is well-known, ${\cal Q}^{+}$ has the following regularity (first found and proved in \cite{Lions1}):
\be\int_{{\bR}}M_1(|\xi|)|\wh{{\cal Q}^{+}(\Psi)}(\xi)|^2{\rm d}\xi
\le (2\pi)^{3}\|\Psi\|_{L^2}^2\lb{Four1}\ee
where the function $M_1(r)\in (0,\infty]$ is defined explicitly through $B_0(z,\sg)$ and
satisfies
$\inf\limits_{r\ge 0}M_1(r)>0, \lim\limits_{r\to\infty}M_1(r)=\infty$;
see Theorem 1 (with $\rho(r)\equiv 1$) in \cite{Lu1998}, see also \cite{BD},\cite{Wen1997}.
By replacing $M_1(r)$ with $\min\{M_1(r), 1+r\}$ we can also assume that
$M_1(r)<\infty$ for all $r\ge 0$.
In application of (\ref{Four1}), it is convenient to use the non-decreasing function
$M_*(R)=\inf\limits_{r\ge R}M_1(r), R>0.$ We have
\bes&&\inf\limits_{R> 0}M_*(R)>0, \quad \lim_{R\to\infty}M_*(R)=\infty,\lb{Decay1}\\
&&\Big(\int_{|\xi|>R}|\wh{{\cal Q}^{+}(\Psi)}(\xi)|^2{\rm d}\xi
\Big)^{1/2}
\le \fr{(2\pi)^{3/2}}{\sqrt{M_*(R)}}
\|\Psi\|_{L^2},\quad R>0.\dnumber\lb{decay2}\ees
If in addition to (\ref{FF}), assume that $F, \wh{F}\in L^1({\bR})$, then
a decay estimate like (\ref{decay2}) holds also for $\wh{{\cal Q}^{+}(\Psi| F)}(\xi)$:
\be\Big(\int_{|\xi|>R}|\wh{{\cal Q}^{+}(\Psi| F)}(\xi)|^2
{\rm d}\xi\Big)^{1/2}
\le {\cal E}_{F}(R)\|\Psi\|_{L^2},\quad R>0\lb{decay3}\ee
where
\be {\cal E}_{F}(R)=\fr{\|\wh{F}\|_{L^1}}{(2\pi)^{3/2}\sqrt{M_*(R/2)}}
+\fr{ \|A_0(|\cdot|)\|_{L^2}}{(2\pi)^{3/2}}\int_{|\eta|> R/2}|\wh{F}(\eta)|
{\rm d}\eta,\quad R>0. \lb{decay4}\ee
To prove this we insert the equality
$F(v_*)=(2\pi)^{-3}\int_{{\bR}}\wh{F}(\eta)e^{{\rm i}\la \eta, v_*\ra}{\rm d}\eta$ (a.e. $v_*\in {\bR}$)
into ${\cal Q}^{+}(\Psi| F)$ to obtain
\beas{\cal Q}^{+}(\Psi| F)(v)=(2\pi)^{-3}\int_{{\bR}}\wh{F}(\eta)e^{-{\rm i}\la \eta, v\ra}
{\cal Q}^{+}(\Psi_{\eta})(v){\rm d}\eta\eeas
where
$\Psi_{\eta}(v,v_*)=\Psi(v,v_*)e^{{\rm i}\la \eta, v+v_*\ra}.$
From this we have
\beas&&\wh{{\cal Q}^{+}(\Psi| F)}(\xi)
=(2\pi)^{-3}\int_{{\bR}}\wh{F}(\eta)
\wh{{\cal Q}^{+}(\Psi_{\eta})}(\xi+\eta){\rm d}\eta,
\\
&& |\wh{{\cal Q}^{+}(\Psi| F)}(\xi)|\le
(2\pi)^{-3}\int_{{\bR}}|\wh{F}(\eta)|
|\wh{{\cal Q}^{+}(\Psi_{\eta})}(\xi+\eta)|{\rm d}\eta.\eeas
Then applying (a derivation of) Minkowski inequality we have, for any $R>0$,
\bes&&\Big(\int_{|\xi|>R}|\wh{{\cal Q}^{+}(\Psi| F)}(\xi)|^2
{\rm d}\xi\Big)^{1/2}
\le (2\pi)^{-3}\int_{{\bR}}|\wh{F}(\eta)|\Big(\int_{|\xi|>R}|\wh{{\cal Q}^{+}(\Psi_{\eta})}(\xi+\eta)|^2{\rm d}\xi\Big)^{1/2}
{\rm d}\eta\nonumber\\
&&=(2\pi)^{-3}\int_{|\eta|\le R/2}|\wh{F}(\eta)|\Big(\int_{|\xi-\eta|>R}|\wh{{\cal Q}^{+}(\Psi_{\eta})}(\xi)|^2{\rm d}\xi\Big)^{1/2}
{\rm d}\eta\nonumber\\
&& +(2\pi)^{-3}\int_{|\eta|> R/2}|\wh{F}(\eta)|\Big(\int_{|\xi-\eta|>R}|\wh{{\cal Q}^{+}(\Psi_{\eta})}(\xi)|^2{\rm d}\xi\Big)^{1/2}
{\rm d}\eta.\lb{QFR}\ees
For the first term in the right hand side of (\ref{QFR}),
since $|\eta|\le R/2, |\xi-\eta|>R$ implies $|\xi|>R/2$, it follows from
(\ref{decay2}) that the inner integral is less than or equal to
$$\Big(\int_{|\xi|>R/2}|\wh{{\cal Q}^{+}(\Psi_{\eta})}(\xi)|^2{\rm d}\xi\Big)^{1/2}
\le \fr{(2\pi)^{3/2}}{\sqrt{M_*(R/2)}}
\|\Psi_{\eta}\|_{L^2}
=\fr{(2\pi)^{3/2}}{\sqrt{M_*(R/2)}}
\|\Psi\|_{L^2}$$
where we used the identity $|\Psi_{\eta}(v,v_*)|\equiv |\Psi(v,v_*)|$.
The estimate for the second term in the right hand side of (\ref{QFR}) follows from
the identity
$\|\wh{f}\|_{L^2}=(2\pi)^{3/2}\|f\|_{L^2}$, (\ref{QF3}), and
$|\Psi_{\eta}|\equiv |\Psi|$. This proves
(\ref{decay3}).

\begin{lemma}\lb{lemma4.3}
Let ${\cal Q}^{+}, {\cal Q}^{+}(\cdot |F)$ be defined in (\ref{calQ}),(\ref{calQF}),
(\ref{FF}),
let $\{\psi_n\}_{n=1}^{\infty}\subset L^{\infty}\cap L^1({\bR})$,
$\{\Psi_n\}_{n=1}^{\infty}\subset L^{\infty}\cap L^1({\bRR})$
satisfy $\sup\limits_{n\ge 1}\{\|\psi_n\|_{L^{\infty}}, \|\psi_n\|_{L^1},
\|\Psi_n\|_{L^{\infty}}, \|\Psi_n\|_{L^1}
\}<\infty$. Then:

 (1) If\, $\psi_n \rightharpoonup 0\,\,\,(n\to\infty)\,$ weakly
 in\,$L^1({\bR})$,
 then
 \be\lim_{n\to\infty}\int_{{\bR}}\psi_n(v){\cal Q}^{+}(\Psi_n)(v){\rm d}v=\lim_{n\to\infty}\int_{{\bR}}\psi_n(v){\cal Q}^{+}(\Psi_n|F)(v){\rm d}v=0.\lb{Lim1}\ee

(2) If\, $\Psi_n \rightharpoonup 0\,\,(n\to\infty)\,$ weakly
 in $L^1({\bRR})$, then
 \be\lim_{n\to\infty}\|{\cal Q}^{+}(\Psi_n)\|_{L^1}=\lim_{n\to\infty}\|{\cal Q}^{+}(\Psi_n|F)\|_{L^1}=0.\lb{Lim2}\ee
\end{lemma}
\begin{proof}  We need only to prove the lemma for ${\cal Q}^{+}(\Psi_n|F)$.
Let $C=\sup\limits_{n\ge 1}\{\|\psi_n\|_{L^{\infty}}, \|\psi_n\|_{L^1},
\|\Psi_n\|_{L^{\infty}},\\  \|\Psi_n\|_{L^1_2}
\}$.
We may assume that $F, \psi_n, \Psi_n$ are real functions.

(1): Suppose
$\psi_n \rightharpoonup 0$ weakly
 in $L^1({\bR})$. Then $\wh{\psi_n}(\xi)\to 0\,\,(n\to\infty)$ for all $\xi\in{\bR}$.
 We first assume in addition that $F,\wh{F}\in L^1({\bR})$. Then
 using Parseval equality we have for any $R>0$
\bes&& (2\pi)^3\Big|\int_{{\bR}}\psi_n(v){\cal Q}^{+}(\Psi_n|F)(v){\rm d}v
\Big|
=\Big|\int_{{\bR}}\overline{\wh{\psi_n}(\xi)}\wh{{\cal Q}^{+}(\Psi_n|F)}(\xi){\rm d}\xi\Big|
\nonumber\\
&&\le\int_{|\xi|\le R}|\wh{\psi_n}(\xi)||\wh{{\cal Q}^{+}(\Psi_n|F)}(\xi)|{\rm d}\xi
+ \int_{|\xi|>R}|\wh{\psi_n}(\xi)||\wh{{\cal Q}^{+}(\Psi_n|F)}(\xi)|{\rm d}\xi\lb{QFW}\ees
and
\beas&&
\sup_{n\ge 1, \xi\in {\bR}}\max\{|\wh{\psi_n}(\xi)|, |\wh{{\cal Q}^{+}(\Psi_n|F)}(\xi)|
\}\le \max\{C, \underline{A}_*C\}.\eeas
By dominated convergence theorem we have
$$\lim_{n\to\infty} \int_{|\xi|\le R}|\wh{\psi_n}(\xi)||\wh{{\cal Q}^{+}(\Psi_n|F)}(\xi)|{\rm d}\xi
=0\quad (\forall\, 0<R<\infty).$$
And using Cauchy-Schwartz inequality and (\ref{decay3}),(\ref{decay4}) gives
\beas\sup_{n\ge 1}
\int_{|\xi|>R}|\wh{\psi_n}(\xi)||\wh{{\cal Q}^{+}(\Psi_n|F)}(\xi)|{\rm d}\xi
\le C^2
{\cal E}_F(R).\eeas
Thus the second limit in (4.20) follows from (\ref{QFW}) by first letting $n\to\infty$ and then letting $R\to\infty$.

Next for general case we use convolution approximation: for any $0<\dt, R<\infty$,
let $F_R(v)=F(v){\bf 1}_{\{|v|\le R\}}, F_{R,\dt}=J_{\dt}*F_R$ where
$J_{\dt}(v)=\dt^{-3}J(\dt^{-1}v), J(v)=(2\pi)^{-3/2}e^{-|v|^2/2},\dt>0$. Then
$F_{R,\dt}, \wh{F_{R,\dt}}\in L^1({\bR})$, $\|F_{R}-F_{R,\dt}\|_{L^1}\to 0\,(\dt\to 0)$, and using the assumption $\|\Psi_n\|_{L^{\infty}},\|\psi_n\|_{L^{\infty}}\le C$ and $\|F\|_{L^{\infty}}\le 1$ we have
\beas&& \Big|\int_{{\bR}}\psi_n(v){\cal Q}^{+}(\Psi_n|F)(v){\rm d}v-\int_{{\bR}}\psi_n(v){\cal Q}^{+}(\Psi_n|F_{R,\dt})(v){\rm d}v
\Big|\\
 &&\le C \int_{{\bRR}}|\psi_n(v)|A_0(|v-v_*|)\big|F(v_*)-F_{R}(v_*)\big|
 {\rm d}v{\rm d}v_*
\\
&&+C\int_{{\bRR}}|\psi_n(v)|A_0(|v-v_*|)
\big|F_R(v_*)-F_{R,\dt}(v_*)\big|
 {\rm d}v{\rm d}v_*\\
 &&
\le \int_{{\bR}}|\psi_n(v)|\Big(\int_{|v_*|>R}A_0(|v-v_*|)
 {\rm d}z\Big){\rm d}v
+
C^2\underline{A}_*\|F_R-F_{R,\dt}\|_{L^1},\eeas
and
\beas&&\int_{{\bR}}|\psi_n(v)|\Big(\int_{|v_*|>R}A_0(|v-v_*|)
 {\rm d}z\Big){\rm d}v=\int_{{\bR}}|\psi_n(v)|\Big(\int_{|v-z|>R}A_0(|z|)
 {\rm d}z\Big){\rm d}v\\
&&\le {\cal E}(R):=C\int_{|z|>R/2}A_0(|z|)
 {\rm d}z
+\|A_0(|\cdot|)\|_{L^1}\sup_{n\ge 1}\int_{|v|> R/2}|\psi_n(v)|{\rm d}v.
\eeas
Thus we obtain
\bes&& \Big|\int_{{\bR}}\psi_n(v){\cal Q}^{+}(\Psi_n|F)(v){\rm d}v\Big|\nonumber\\
&&\le {\cal E}(R)
+C^2\underline{A}_*\|F_R-F_{R,\dt}\|_{L^1}+
\Big|\int_{{\bR}}\psi_n(v){\cal Q}^{+}(\Psi_n|F_{R,\dt})(v){\rm d}v
\Big|.\lb{calQPF} \ees
By $F_{R,\dt}, \wh{F_{R,\dt}}\in L^1({\bR})$ we have
$$\lim_{n\to\infty}
\Big|\int_{{\bR}}\psi_n(v){\cal Q}^{+}(\Psi_n|F_{R,\dt})(v){\rm d}v
\Big|=0.$$
Since $\{\psi_n\}_{n=1}^{\infty}$ is weakly compact in $L^1({\bR})$, it follows from
Dunford-Pettis criterion of $L^1$-relatively weak compactness (\cite{DS}) that
$\sup\limits_{n\ge 1}\int_{|v|> R/2}|\psi_n(v)|{\rm d}v\to 0\,(R\to\infty)$.
Thus ${\cal E}(R)\to 0$ as $R\to\infty$.
By first taking upper limit $\limsup\limits_{n\to\infty}$ to the left hand side of
 (\ref{calQPF}), and then letting $\dt\to 0$ and finally letting $R\to\infty$ we deduce
the second limit in (\ref{Lim1}).

(2): Suppose $\Psi_n \rightharpoonup 0$ weakly
 in $L^1({\bRR})$. Then ${\cal Q}^{+}(\Psi_n|F) \rightharpoonup 0$ weakly
 in $L^1({\bR})$. In fact for any $\psi\in L^{\infty}({\bR})$,
$$\int_{{\bR}}\psi {\cal Q}^{+}(\Psi_n|F)
{\rm d}v=\int_{{\bRR}}\Psi_n(v,v_*)\vp(v,v_*){\rm d}v{\rm d}v_*\to 0\,\,\,
(n\to\infty)$$
where
$\vp(v,v_*)=\int_{{\bS}}B_0(v-v_*,\sg)F(v_*')\psi(v'){\rm d}\sg$
which is obviously in $L^{\infty}({\bRR})$.

Since $\sup\limits_{n\ge 1}\{\|{\cal Q}^{+}(\Psi_n|F)\|_{L^{\infty}},
\|{\cal Q}^{+}(\Psi_n|F)\|_{L^1}\}\le \max\{\|A_0(|\cdot|)\|_{L^1}C, \,
\underline{A}_*C\}<\infty$, it follows from part (1) with
$\psi_n ={\cal Q}^{+}(\Psi_n|F)$ that  $\|{\cal Q}^{+}(\Psi_n|F)\|_{L^2}
 \to 0\,(n\to\infty)$.

 Now again for any $0<R<\infty$ we have
\beas
\|{\cal Q}^{+}(\Psi_n|F)\|_{L^1}
\le (4\pi R^3/3)^{1/2}\|{\cal Q}^{+}(\Psi_n|F)\|_{L^2}+\int_{|v|>R}{\cal Q}^{+}(\Psi_n|F)(v)|{\rm d}v.
\eeas
Since $\{{\cal Q}^{+}(\Psi_n|F)\}_{n=1}^{\infty}$ is $L^1$-weakly compact, it follows from
criterion of $L^1$-relatively weak compactness that
$\sup\limits_{n\ge 1}\int_{|v|>R}{\cal Q}^{+}(\Psi_n|F)(v)|{\rm d}v\to 0$ as $R\to\infty$.
This proves the second limit in (\ref{Lim2}).
\end{proof}
\vskip2mm

The following lemma presents some basic convergence properties
that are used in our proof of convergence to equilibrium.

\begin{lemma}\lb{lemma4.4} (1)
Let $\{f_n\}_{n=1}^{\infty}\subset L^1({\bR})$ be relatively weakly compact in
$L^1({\bR})$ and  let $\{g_n\}_{n=1}^{\infty}\subset L^{\infty}({\bR})$
satisfy $\sup\limits_{n\ge 1}\|g_n\|_{L^{\infty}}<\infty$ and
$\lim\limits_{n\to\infty}g_n(v)=0$ for a.e. $v\in {\bR}$.
Then $\lim\limits_{n\to\infty}\|f_ng_n\|_{L^1}=0$.

(2) Let $\{f_n\}_{n=1}^{\infty}\subset L^1({\bR})$ be relatively weakly compact in
$L^1({\bR})$ and let $\{g_n\}_{n=1}^{\infty}\subset L^1({\bR})$
satisfy $g_n\rightharpoonup 0\,(n\to\infty)$ weakly in $L^1({\bR})$.
Then $f_n\otimes g_n\rightharpoonup 0, g_n\otimes f_n\rightharpoonup 0\,(n\to\infty)$  weakly in $L^1({\bRR})$.

(3) $\{f_n\}_{n=1}^{\infty}, \{g_n\}_{n=1}^{\infty}\subset L^1({\bR})$ and $f,g\in L^1({\bR})$
be such that  $f_n\rightharpoonup f, g_n\rightharpoonup g\,(n\to\infty)$  weakly in $L^1({\bR})$.
Then $f_n\otimes g_n\rightharpoonup f\otimes g\,(n\to\infty)$  weakly in $L^1({\bRR})$.
\end{lemma}
\begin{proof} Part (1) follows from the criterion of
$L^1$-relatively weak compactness (see e.g. \cite{DS}) and dominated convergence theorem.
 Part (2) is an application of part (1). Part (3) follows from part (2).
 \end{proof}
\vskip2mm

We also need strong or weak Lipschitz continuity in time for mild or weak solutions:

\begin{lemma}\lb{lemma4.2} (1)
Let $B(z,\sg)$ satisfy (A1),(A2), let $f_0\in L^1_2({\bR})$ satisfy
$0\le f_0\le 1$ on ${\bR}$.
Let $f$ be the unique conservative mild solution of Eq.(BFD) with
the initial datum $f_0$. Then
\be\|f(t_1)-f(t_2)\|_{L^1},
 \|Q_0(f)(t_1)-Q_0(f)(t_2)\|_{L^1}\le C|t_1-t_2|\quad \forall\,t_1, t_2\in [0,\infty)\lb{Diff1}
 \ee
 where $Q_0(\cdot)$ is the collision operator corresponding to $B_0(z,\sg)$
 defined  (\ref{smaller}) for $0\le \gm\le 1$, $0<C<\infty$
 depends only on $\int_{0}^{\pi}b^*(\cos\theta)\sin(\theta){\rm d}\theta,
 \int_{0}^{\pi}\underline{b}_*(\cos\theta)\sin(\theta){\rm d}\theta $ and $\|f_0\|_{L^1_2}$.

 (2)
Let $B(z,\sg)$ satisfy (A1),(A3), let $f_0\in L^1_2({\bR})$ satisfy
$0\le f_0\le 1$ on ${\bR}$.
Let $f$ be a conservative weak solution of Eq.(BFD) with
the initial datum $f_0$. Then
\be\Big|\int_{{\bR}}\big(f(t_1,v)-f(t_2,v)\big)\vp(v)
{\rm d}v\Big|\le C_{\vp}|t_1- t_2|\quad \forall\, t_1, t_2\in[0,\infty),\,\,
\forall\,\vp\in C^2_b({\bR})\ee
where $C_{\vp}$ depends only on $\|\p^2\vp\|_{L^{\infty}}, \int_{0}^{\pi}b^*(\cos\theta)
\sin^2(\theta){\rm d}\theta, \gm, \|f_0\|_{L^1}, \|f_0\|_{L^1_2}$.
 \end{lemma}

\begin{proof} (1): Since $0\le \gm\le 1$, the
inequality $\|f(t_1)-f(t_2)\|_{L^1}\le C|t_1-t_2|$ follows from
$|v-v_*|^{\gm}\le \la v\ra ^{\gm}+\la v_*\ra^{\gm}$ and the conservation of
mass and energy: $\|f(t)\|_{L^1_{\gm}}\le \|f_0\|_{L^1_2}.$
To prove the second inequality in (\ref{Diff1}), we temporarily denote
$f(v)=f(t_1,v), g(v)=f(t_2,v)$. Then
\bes&&|\Pi_{{\rm F}}(f)-\Pi_{{\rm F}}(g)|
\le f'f_*'(|f-g|+|f_*-g_*|)+|f'-g'|f_*'+g'|f_*'-g_*'|\nonumber\\
&&+ff_*(|f'-g'|+|f_*'-g_*'|)+|f-g|f_*+g|f_*-g_*|
\Longrightarrow\nonumber\\
&&\|Q_0(f)(t_1)-Q_0(f)(t_2)\|_{L^1}
\le \int_{{\bRRS}}B_0(v-v_*,\sg)
|\Pi_{{\rm F}}(f)-\Pi_{{\rm F}}(g)|{\rm d}\sg{\rm d}v_*{\rm d}v
\nonumber \\
&&\le 4\int_{{\bR}}|f-g|\Big(\int_{{\bRS}}
B_0(v-v_*,\sg) f'f_*'{\rm d}\sg{\rm d}v_*\Big){\rm d}v
\nonumber\\
&&+2\underline{A}_*\int_{{\bRR}}\big(
|f-g|f_*+g|f_*-g_*|\big){\rm d}v_*{\rm d}v.\lb{Q0}
\ees
For the first integral in the right hand side of (\ref{Q0}),
consider $f'f_*'\le f'{\bf 1}_{\{\pi/2< \theta\le\pi\}}
+ f_*'{\bf 1}_{\{0\le \theta\le\pi/2\}}$ and
use Lemma \ref{lemmA3.1} and notice that $B_0(z,\sg)\le \underline{b}_*(\cos\theta)$
to deduce
\beas&& \int_{{\bRS}}B_0(v-v_*,\sg) f'f_*'{\rm d}\sg{\rm d}v_*
\\
&&\le 2\pi\int_{\bR}
 f(v_*)\Bigg\{\int_{\pi/2}^{\pi}\underline{b}_*(\cos\theta)\frac{\sin(\theta)}{\sin^{3}(\theta/2)}
 {\rm d}\theta+
\int_{0}^{\pi/2}\underline{b}_*(\cos\theta)\frac{\sin(\theta)}{\cos^3(\theta/2)}{\rm d}\theta\Bigg\}{\rm d}v_*
\\
&&\le C_0\underline{A}_*\int_{\bR}
 f(v_*){\rm d}v_*=C_0\underline{A}_*\|f_0\|_{L^1}.\eeas
 Thus we obtain
\beas&&\|Q_0(f)(t_1)-Q_0(f)(t_2)\|_{L^1}
\le C_0\underline{A}_*\int_{{\bRR}}\big(
|f-g|f_*+g|f_*-g_*|\big){\rm d}v_*{\rm d}v.
\\
&&\le C_0\underline{A}_*\|f_0\|_{L^1}\|f(t_1)-f(t_2)\|_{L^1}
\le C |t_1-t_2|.\eeas

(2): Recalling definition of weak solutions we have
$$\int_{{\bR}}\big(f(t_1,v)-f(t_2,v)\big)\vp(v)
{\rm d}v=\int_{t_2}^{t_1}{\rm d}t\int_{{\bR}}
B\Pi_{{\rm F}}(f)\Dt\vp{\rm d}\sg{\rm d}v_*{\rm d}v$$
and so using (\ref{qq3}) we obtain
\beas&&\Big|\int_{{\bR}}\big(f(t_1)-f(t_2)\big)\vp
{\rm d}v\Big|\le A^*\|\p^2\vp\|_{L^{\infty}}
\int_{t_2}^{t_1}{\rm d}t\int_{{\bRR}}ff_*|v-v_*|^{2-|\gm|}{\rm d}v_*{\rm d}v
\\
&&\le C_{\vp}|t_1-t_2|.\eeas
\end{proof}

{\bf Proof of Theorem \ref{theorem1.1} and Theorem \ref{theorem1.2}.}	 Let $B, f_0$ and
the (mild or weak) conservative solution $f$ of Eq.(BFD) be given in the two theorems respectively.
For  Theorem \ref{theorem1.1}, it has been proven in \cite{Lu2001},\cite{LW2003} that
the conservative mild solution $f$ exists uniquely and satisfies
the entropy identity (\ref{entropy-identity}) and the moment estimates: if $0<\gm\le 1$, then for any $s>2$
\be \sup_{t\ge 1}\|f(t)\|_{L^1_s}<\infty\lb{moment4.1}\ee
 while if $\gm=0$, then (\ref{moment4.1}) also holds for $s>2$ provided that $f_0\in L^1_s({\bR})$.
For Theorem \ref{theorem1.2}, the existence of the weak solution $f$ satisfying the entropy inequalities (\ref{weak-entropy1})-(\ref{weak-entropy4}) and the moment estimates (\ref{a-bounds}) has been proven in Section 3. Thus in the following we need only to prove the
convergence part. Let
$F$ be the equilibrium having the same mass, mean velocity, and energy as $f_0$.
Let $M_0, M_2, v_0$ be defined in (\ref{mac}).
If $M_2(M_0)^{-5/3}=(3/5)(3/(4\pi))^{2/3}$, then in view of
(\ref{equilib2}) and the conservation of mass, momentum, and energy we conclude that
$f(t,v)=F(v)={\bf 1}_{\{|v-v_0|\le R\}}$ for all $t\in [0,\infty)$ and for all $v\in {\bR}\setminus Z_t$ where
$Z_t$ is null set for every $t\ge 0$. So in this case we have
$\|f(t)-F\|_{L^1_2}\equiv 0, t\in [0,\infty)$. Thus
in the following we assume that $M_2(M_0)^{-5/3}>(3/5)(3/(4\pi))^{2/3}$
and we need only to prove the strong and averaged strong convergence of $f(t,\cdot)$ to $F$
as stated in the two theorems. Note that in this case, the equilibrium $F$ given by (\ref{equilib1}) satisfies $\wh{F}\in L^1({\bR})$.

{\bf Step1.} Let $Q_0(f)=\int_{{\bRS}}B_0\Pi_{{\rm F}}(f){\rm d}\sg{\rm d}v_*$ be the collision operator (\ref{Q})
 corresponding to $B_0(z,\sg)$ and let
\beas {\cal N}(f)(v)=\int_{{\bRS}}B_0\big(f'f_*'+f_*(1-f'-f_*')\big){\rm d}\sg{\rm d}v_*.\eeas
Then
$$Q_0(f)= {\cal Q}^{+}(f\otimes f|1-f)-f {\cal N}(f).$$
Since $F$ is an equilibrium,  we have
$Q_0(F)=0$,
i.e. $F{\cal N}(F)={\cal Q}^{+}(F\otimes F|1-F).$
So
\beas (f-F){\cal N}(F)
=f( {\cal N}(F)-{\cal N}(f))+{\cal Q}^{+}(f\otimes f|1-f)-
{\cal Q}^{+}(F\otimes F|1-F)-Q_0(f).\eeas
Let us write
$$f=f(t,\cdot),\quad h=f-F,\quad g=fh.$$
Then $f=F+h$ and
\beas &&(f-F)^2{\cal N}(F)
=g\big( {\cal N}(F)-{\cal N}(F+h)\big)\\
&&+h\big({\cal Q}^{+}\big((F+h)\otimes (F+h)|1-F-h)-{\cal Q}^{+}(F\otimes F|1-F)\big)-hQ_0(f)\eeas
Further calculation deduce
\be\int_{{\bR}}(f(t)-F)^2{\cal N}(F){\rm d}v=
-\int_{{\bR}}h{\cal Q}^{+}(h\otimes h|h){\rm d}v
-\int_{{\bR}}hQ_{0}(f){\rm d}v+W(t)\lb{NW1}\ee
where
\beas
W(t)&=&-\int_{{\bR}}g{\cal Q}^{+}((F+f)\otimes h){\rm d}v
-\int_{{\bRR}}(g\otimes h)A_0(|v-v_*|){\rm d}v_*{\rm d}v\\
&+&
2\int_{{\bR}}f{\cal Q}^{+}(g\otimes h){\rm d}v
-\int_{{\bR}}f{\cal Q}^{+}(h\otimes h|F){\rm d}v
\\
&+&2\int_{{\bR}}h\Big\{{\cal Q}^{+}(g\otimes F)+ {\cal Q}^{+}(F\otimes h)
-{\cal Q}^{+}(f\otimes h|F)+\fr{1}{2}{\cal Q}^{+}(h\otimes h)\Big\}
{\rm d}v
.\eeas
By Proposition \ref{prop4.1} we have
$$\int_{{\bR}}h{\cal Q}^{+}(h\otimes h|h){\rm d}v\ge 0$$
hence
\be\int_{{\bR}}(f(t,v)-F(v))^2 {\cal N}(F)(v){\rm d}v\le -\int_{{\bR}}hQ_{0}(f){\rm d}v
+W(t),\quad t\ge 0.\lb{NW2}\ee

\noindent {\bf Step2.}
Recalling the equilibrium $F(v)=ae^{-b|v-v_0|^2}/(1+ae^{-b|v-v_0|^2})$ we have
\beas&& F'F_*'+F_*(1-F'-F_*')=\fr{1}{F} F'F_*'(1-F_*)\ge\fr{1}{1+a}F_*,\\
&&
\Longrightarrow\,{\cal N}(F)(v)\ge \fr{1}{1+a}\int_{{\bR}}A_0(|v-v_*|)F(v_*){\rm d}v_*
\\
&&\Longrightarrow\,N_R:=\min\limits_{|v|\le R}{\cal N}(F)(v)>0\quad \forall\, 0<R<\infty.\eeas
We then prove that in order to prove the strong convergence (\ref{strong}) (for $0\le \gm\le 1$)
and the averaged  strong convergence (\ref{averaging strong}) (for $-4\le \gm<0$), it needs only to prove that
\be \lim_{t\to\infty}\int_{{\bR}}(f(t)-F)^2 {\cal N}(F){\rm d}v=0,\quad
\lim_{t\to\infty}\fr{1}{t}\int_{0}^{t}{\rm d}\tau\int_{{\bR}}(f(\tau)-F)^2 {\cal N}(F){\rm d}v=0
\lb{SA}\ee respectively. In fact by identities
$\|f(t)\|_{L^1_2}=\|F\|_{L^1_2}$ and $|f(t,v)-F(v)|=f(t,v)-F(v)+2(F(v)-f(t,v))_{+}$ we have
for any $1\le R<\infty$
\beas&&
\|f(t)-F\|_{L^1_2}=2\int_{{\bR}}(F-f(t))_{+}\la v\ra^2{\rm d}v
\\
&&=2\int_{|v|\le R}(F-f(t))_{+}\la v\ra^2{\rm d}v
+2\int_{|v|>R}(F-f(t))_{+}\la v\ra^2{\rm d}v
\\
&&\le 8\big(\fr{\pi}{3N_R}R^7\big)^{1/2}\Big(\int_{|v|\le R}(F-f(t))^2{\cal N}(F){\rm d}v
\Big)^{1/2}+\fr{2}{R}\|F\|_{L^1_3}.\eeas
This implies the above reduction.

{\bf Step3.}
Using (\ref{pi6}) and Cauchy-Schwartz inequality we have
\be
\int_{{\bRRS}}B_0(v-v_*,\sg)|\Pi_{{\rm F}}(f)(t,v,v_*,\sg)|
{\rm d}\sg{\rm d}v_*{\rm d}v
\le C_0\|f_0\|_{L^1}\sqrt{D(f,q)(t)}
\lb{control1}\ee
and in particular
\be\|Q_0(f)(t)\|_{L^1}\le C_0\|f_0\|_{L^1}\sqrt{D(f,q)(t)},\quad t\ge 0.
\lb{controlq}\ee
Recall that for the case $0\le \gm\le 1$ we have $q=f'f_*'ff_*$ so that
$D(f,q)(t)=D(f)(t)$.

By entropy inequalities (\ref{weak-entropy1}),(\ref{SLeq}) we have
\be\int_{0}^{\infty} D(f,q)(t){\rm d}t\le \sup_{t\ge 0}\big(S(f(t))-S(f_0)\big)<\infty.\lb{entropy3}\ee
We now prove that for any sequence $\{t_n\}_{n=1}^{\infty}\subset [1,\infty)$ satisfying
\be \lim_{n\to\infty}t_n=\infty\quad {\rm and}\quad \sup_{n\ge 1}\|f(t_n)\|_{L^1_s}<\infty\quad{\rm for\,\,some}\,\,s>2 \lb{M3}\ee
there exists a subsequence of $\{t_n\}_{n=1}^{\infty}$, still denote it as  $\{t_n\}_{n=1}^{\infty}$,
such that
$$\lim_{n\to\infty}|W(t_n)|=0$$
and for the case $0\le \gm\le 1$, it also holds
$$\lim_{n\to\infty}\|Q_0(f)(t_n)\|_{L^1}=0.$$
To do this we use entropy dissipation inequalities. Let $\dt_n=\big(\int_{t_n}^\infty D(f,q)(t){\rm d}t+1/n\big)^{1/2}$. Then
$0<\dt_n\to 0\,(n\to\infty)$ and
$
\frac{1}{\dt_n}\int_{t_n}^{t_n+\dt_n}D(f,q)(t){\rm d}t<\dt_n$
so that there exist $\tau_n\in[t_n,t_n+\dt_n]$ such that
\be D(f,q)(\tau_n)\le\dt_n\to 0\quad (n\to\infty).\lb{Dq}\ee
By $0\le f\le 1$ on $[0,\infty)\times {\bR}$ and $\|f(t)\|_{L^1_2}
\equiv\|f_0\|_{L^1_2}$,
there exist a subsequence of $\{(t_n,\tau_n)\}_{n=1}^{\infty}$, still denote it as  $\{(t_n,\tau_n)\}_{n=1}^{\infty}$, and $f_{\infty}, \tilde{f}_{\infty} \in L^1({\bR})$ such that
$f(t_n, \cdot)\rightharpoonup f_\infty, f(\tau_n, \cdot)\rightharpoonup \tilde{f}_\infty$
weakly in $L^1(\bR).$  Since $0\le f\le 1$ on $[0,\infty)\times{\bR}$,
we can modify $f_{\infty}, \tilde{f}_{\infty}$ on
a null set such that $0\le f_{\infty}, \tilde{f}_{\infty}\le 1$ on ${\bR}.$ Also by the
second condition in (\ref{M3}) we have
\be\int_{\bR}f_\infty(v)(1, v, |v|^2){\rm d}v=\int_{\bR}f_0(v)(1, v, |v|^2){\rm d}v.\lb{MME}\ee
Using Lemma \ref{lemma4.2} we have for the case $0\le \gm\le 1$ that $f(t_n,\cdot)-f(\tau_n,\cdot)\to 0\, (n\to\infty)$ in $L^1$ and so $f_{\infty}=\tilde{f}_{\infty}$ a.e. on ${\bR}$; while for the case $-4\le \gm<0$ we have, as $n\to\infty$,
\beas&&\Big|\int_{{\bR}}(f(t_n,v)-f(\tau_n,v))\vp(v)
{\rm d}v\Big|\le C_{\vp}|t_n-\tau_n|\le C_{\vp}\dt_n\to 0\quad
\forall\,\vp\in C^2_b({\bR})\\
&& \Longrightarrow\,\int_{{\bR}}(f_{\infty}-\tilde{f}_{\infty})\vp{\rm d}v=0\qquad
\forall\,\vp\in C^2_b({\bR})\eeas
and so we also conclude $\wt{f}_{\infty}=f_{\infty}$ a.e. on ${\bR}$.

For further estimate we write $f_n(v)=f(\tau_n,v), h_n(v)=f(\tau_n,v)-f_{\infty}(v)$.
By convexity $y_0^2\le y^2-2y_0(y-y_0)$ we have
\bes&&\int_{{\bRRS}} B_0 \big(\Pi_{{\rm F}}(f_{\infty})\big)^2
{\rm d}\sg{\rm d}v{\rm d}v_*
\le \int_{{\bRRS}} B_0 \big(\Pi_{{\rm F}}(f_n)\big)^2
{\rm d}\sg{\rm d}v{\rm d}v_*
\nonumber\\
&&-2 \int_{{\bRRS}} B_0\Pi_{{\rm F}}(f_\infty)
\big(\Pi_{{\rm F}}(f_n)-\Pi_{{\rm F}}(f_\infty)\big)
{\rm d}\sg{\rm d}v{\rm d}v_*.
\lb{FE}\ees
Using $(\Pi_{{\rm F}}(f_n))^2
\le |\Pi_{{\rm F}}(f_n)|$, (\ref{control1}) and (\ref{Dq}) we have
\be\int_{{\bRRS}} B_0 \big(\Pi_{{\rm F}}(f_n)\big)^2
{\rm d}\sg{\rm d}v{\rm d}v_*
\le C_0\|f_0\|_{L^1}\sqrt{D(f,q)(\tau_n)}\to 0\quad (n\to\infty).\lb{4.37}\ee
Observe that the quartic terms $ff_*f'f_*'$ in
$\Pi_{{\rm F}}(f)$ cancel each other:
$\Pi_{{\rm F}}(f)=f'f_*'(1-f-f_*)-ff_*(1-f'-f_*').$
Using this we compute
\beas&& -\int_{{\bRRS}} B_0\Pi_{{\rm F}}(f_\infty)
\big(\Pi_{{\rm F}}(f_n)-\Pi_{{\rm F}}(f_\infty)\big)
{\rm d}\sg{\rm d}v{\rm d}v_*
\\
&&=\int_{{\bRR}}
\big(f_n\otimes {f_n}-f_{\infty}\otimes {f_\infty}
\big)\vp(v,v_*){\rm d}v{\rm d}v_*
\\&&+4\int_{{\bR}}h_n (1-f_{\infty})
{\cal Q}^{+}\big( (f_{\infty}f_n)\otimes
(f_{\infty}f_n)\big){\rm d}v
-4\int_{{\bR}}h_n{\cal Q}^{+}\big( (f_{\infty}f_n)\otimes
(f_{\infty}f_n)|f_{\infty}\big){\rm d}v\\
&&-4\int_{{\bR}}h_n f_{\infty}{\cal Q}^{+}(f_n\otimes f_n| f_{\infty})
{\rm d}v+
8\int_{{\bR}}h_nf_{\infty}{\cal Q}^{+}\big( (f_nf_{\infty})\otimes
f_n|f_{\infty}\big){\rm d}v
\eeas
where
$$\vp(v,v_*)=-2\big(1-2f_{\infty}(v))\int_{{\bS}}B_0(v-v_*,\sg)\Pi_{{\rm F}}(f_{\infty})(v,v_*,\sg)
{\rm d}\sg$$
which is in $L^{\infty}({\bRR})$.
By Lemma \ref{lemma4.4} and Lemma \ref{lemma4.3} we
have
$$-\int_{{\bRRS}} B_0\Pi_{{\rm F}}(f_\infty)
\big(\Pi_{{\rm F}}(f_n)-\Pi_{{\rm F}}(f_\infty)\big)
{\rm d}\sg{\rm d}v{\rm d}v_*\to 0\quad (n\to\infty)$$
and so letting $n\to\infty$ we conclude from (\ref{FE}) and (\ref{4.37}) that
$$\int_{{\bRRS}} B_0 \big(\Pi_{{\rm F}}(f_{\infty})\big)^2
{\rm d}\sg{\rm d}v{\rm d}v_*=0.$$
Since $B_0(v-v_*,\sg)>0$ a.e. $(v,v_*,\sg)\in {\bRRS}$, this implies that
$\Pi_{{\rm F}}(f_{\infty})=0$ a.e. $(v,v_*,\sg)\in {\bRRS}$. So
$f_{\infty}$ is an equilibrium of Eq.(BFD). Then together with (\ref{MME}) we conclude
that $f_{\infty}=F$ a.e. in ${\bR}$.
 Thus we conclude $f(t_n,\cdot)\rightharpoonup F.$
Let $f_n(v)=f(t_n,v), h_n(v)=f(t_n,v)-F(v), g_n(v)=f_n(v)h_n(v)$.
Then, by definition of $W(t)$ defined above, we have
\beas
W(t_n)&=&-\int_{{\bR}}g_n{\cal Q}^{+}((F+f_n)\otimes h_n){\rm d}v
-\int_{{\bRR}}(g_n\otimes h_n)A_0(|v-v_*|){\rm d}v_*{\rm d}v\\
&+&
2\int_{{\bR}}f_n{\cal Q}^{+}(g_n\otimes h_n){\rm d}v
-\int_{{\bR}}f_n{\cal Q}^{+}(h_n\otimes h_n|F){\rm d}v
\\
&+&2\int_{{\bR}}h_n\Big\{{\cal Q}^{+}(g_n\otimes F)+ {\cal Q}^{+}(F\otimes h_n)
-{\cal Q}^{+}(f_n\otimes h_n|F)+\fr{1}{2}{\cal Q}^{+}(h_n\otimes h_n)\Big\}
{\rm d}v
.\eeas
Since $h_n=f_n-F\rightharpoonup 0$ weakly in $L^1({\bR})$ and (by Lemma \ref{lemma4.4})
$(F+f_n)\otimes h_n\rightharpoonup 0, g_n\otimes h_n\rightharpoonup 0, h_n\otimes h_n\rightharpoonup 0 $  weakly  in $L^1({\bRR})$,
and since all sequences of these functions inside ${\cal Q}^{+}(\cdot), {\cal Q}^{+}(\cdot |F)$
 are bounded both in $L^{\infty}({\bRR})$ and $L^1({\bRR})$,
it follows from Lemma \ref{lemma4.3} that
$W(t_n)\to 0\,(n\to\infty)$. Finally
if $0\le \gm\le 1$, then  using Lemma \ref{lemma4.2},
(\ref{controlq}) and (\ref{Dq}) we have
$$\|Q_0(f)(t_n)\|_{L^1}\le \|Q_0(f)(\tau_n)\|_{L^1}+C|t_n-\tau_n|\le
C\sqrt{D(f,q)(\tau_n)}+
C\dt_n \to 0\quad (n\to\infty).$$
{\bf Step4.} We prove the two convergences in (\ref{SA}) for $0\le \gm\le 1$ and $-4\le \gm<0$
respectively.

\noindent {\bf Case1:} $0\le \gm\le 1.$
In this case the moment estimate (\ref{moment4.1}) holds for some $s>2$. Choose a sequence
$\{t_n\}_{n=1}^{\infty}\subset [1,\infty)$ satisfying $t_n\to\infty\,(n\to\infty)$ such that
$$\lim_{n\to\infty}\big(\|Q_{0}(f)(t_n)\|_{L^1}+|W(t_n)|\big)=
\limsup_{t\to\infty}\big(\|Q_{0}(f)(t)\|_{L^1}+|W(t)|\big).$$
Then using {\bf Step3} there is a subsequence of $\{t_n\}_{n=1}^{\infty}$,
still denote is as $\{t_n\}_{n=1}^{\infty}$, such that
$\lim\limits_{n\to\infty}\big(|Q_{0}(f)(t_n)\|_{L^1}+|W(t_n)|\big)=0$.
Thus
$$\limsup_{t\to\infty}\int_{{\bR}}(f(t)-F)^2 {\cal N}(F){\rm d}v\le \limsup_{t\to\infty}\big(\|Q_{0}(f)(t)\|_{L^1}+|W(t)|\big)=0.$$

\noindent {\bf Case2:} $-4\le\gm<0.$ From (\ref{NW2}) we have for all $t>0$
\be\fr{1}{t}\int_{0}^{t}{\rm d}\tau\int_{{\bR}}(f(\tau)-F)^2 {\cal N}(F){\rm d}v\le
\fr{1}{t}\int_{0}^{t}\|Q_{0}(f)(\tau)\|_{L^1}{\rm d}\tau+\fr{1}{t}\int_{0}^{t}|W(\tau)|{\rm d}\tau.\lb{NW3}\ee
We prove that each term in the right hand side of (\ref{NW3}) tends to zero as $t\to\infty$.
First using (\ref{controlq}) and (\ref{entropy3}) we have
$$\fr{1}{t}\int_{0}^{t}\|Q_0(f)(\tau)\|_{L^1}{\rm d}\tau
\le C_0\|f_0\|_{L^1}\int_{0}^{t}\sqrt{D(f,q)(\tau)}{\rm d}\tau
\le C\fr{1}{\sqrt{t}}\to 0\,\,\, (t\to\infty).$$
Next, since $C_1:=\sup\limits_{t\ge 0}|W(t)|<\infty$,
to prove $\frac{1}{t}\int_0^t |W(\tau)|{\rm d}\tau\to 0
\,(t\to\infty)$
it needs only to prove
\be \lim_{t\to\infty}\fr{1}{|I(t)|}\int_{I(t)}|W(\tau)|{\rm d}\tau=0,\quad {\rm where}
\quad I(t)=[\sqrt{t},\, t],\, t\ge 4.\lb{WA}\ee
By averaging estimate for moments we have
$$C_2=\sup_{t\ge 4}\fr{1}{|I(t)|}\int_{I(t)}\|f(\tau)\|_{L^1_{\td{s}}}{\rm d}\tau<\infty,\quad \td{s}=s-|\gm|>2 .$$
For any $\vep>0$, choose a sequence $\{t_n^*\}_{n=1}^{\infty}=\{t_{n,\vep}^*\}_{n=1}^{\infty}
\subset [4,\infty)$ satisfying $t_n^*\to\infty(n\to\infty)$ such that
$$\lim_{n\to\infty}\fr{1}{I(t_n^*)}\int_{I(t_n^*)}
\big(|W(\tau)|-\vep\|f(\tau)\|_{L^1_{\td{s}}}\big){\rm d}\tau
=\limsup_{t\to\infty}\fr{1}{|I(t)|}\int_{I(t)}
\big(|W(\tau)|-\vep\|f(\tau)\|_{L^1_{\td{s}}}\big){\rm d}\tau.$$
For every $n\in\mathbb{N}$ there exists $t_n=t_{n,\vep}\in I(t_n^*)$ such that
\beas\fr{1}{|I(t_n^*)|}\int_{I(t_n^*)}
\big(|W(\tau)|-\vep\|f(\tau)\|_{L^1_{\td{s}}}\big){\rm d}\tau\le |W(t_n)|-\vep\|f(t_n)\|_{L^1_{\td{s}}}.
\eeas
This implies that $\{\|f(t_n)\|_{L^1_{\td{s}}}\}_{n=1}^{\infty}$ is bounded:
\beas
&&\|f(t_n)\|_{L^1_{\td{s}}}\le
\frac{1}{\vep}|W(t_n)|+\fr{1}{|I(t_n^*)|}\int_{I(t_n^*)}
\|f(\tau)\|_{L^1_{\td{s}}}{\rm d}\tau
\le \fr{C_1}{\vep}+C_2,\,\, n=1,2,3,...\,.\eeas
Thus by {\bf Step3} there exists a subsequence of $\{t_n\}_{n=1}^{\infty}$, still denote it as $\{t_n\}_{n=1}^{\infty},$ such that
$|W(t_n)|\to 0$ as $n\to\infty$, and we then deduce
\beas
&&\limsup_{t\to\infty}\frac{1}{|I(t)|}\int_{I(t)}|W(\tau)|{\rm d}\tau\\
&&\le \limsup_{t\to\infty}\frac{1}{|I(t)|}
\int_{I(t)}\Big(|W(\tau)|-\vep\|f(\tau)\|_{L^1_{\td{s}}}\Big){\rm d}\tau+
\vep \limsup_{t\to\infty}\frac{1}{|I(t)|}
\int_{I(t)}\|f(\tau)\|_{L^1_{\td{s}}}{\rm d}\tau\\
&&\le \lim_{n\to\infty}|W(t_n)|+\vep C_2=\vep C_2.
\eeas
Letting $\vep\to 0$ we obtain (\ref{WA}).
This finishes the proof of Theorem \ref{theorem1.1} and Theorem \ref{theorem1.2}.
$\hfill\Box$
\\

\begin{center}\section {Appendix}\end{center}

Here we show some collision kernels $B(z,\sg)$ of the form (\ref{weak-coupling}) that satisfy the assumptions \{(A1),(A2),(A4)\} and \{(A1),(A3),(A4)\} respectively. A counterexample for the
nonnegativity (\ref{positivity}) is also given.

{\bf 5.1.}
Let ${\cal S}({\bR})$ be the class of Schwartz functions on ${\bR}$.
The (generalized) Fourier transform
$\wh{\phi}(|\xi|)$ of a spherically symmetric
 interaction potential $\phi(|x|)$ is defined by
$$
 \int_{{\bR}}\wh{\phi}(|\xi|)\psi(\xi){\rm d}\xi=\int_{{\bR}}\phi (|x|)\wh{\psi}(x){\rm d}{\bf x}\qquad \forall\, \psi\in {\cal S}({\bR})$$
where
$\wh{\psi}(x)=
\int_{{\bR}}\psi(\xi)e^{-{\rm i}\la x,\xi\ra}{\rm d}\xi.$

For the inverse power-law potential $\phi(|x|)=|x|^{-\alpha}$ with $0<\alpha<3$,
its generalized Fourier transform is
$$\wh{\phi}(|\xi|)=\fr{2^{3-\alpha}\pi^{\fr{3}{2}}\Gm(\fr{3-\alpha}{2})}{\Gm(\fr{\alpha}{2})}
|\xi|^{-3+\alpha},\quad \xi\in{\bR}\setminus\{0\}.$$
See e.g. Lemma 1 in \S 1.1 of Chapter 5 and Lemma in  \S 3.3  of Chapter 3 in \cite{Stein}
and so from  (\ref{weak-coupling}) we have
\bes&& B(z,\sg)=|z|^{\gamma }b(\cos\theta),\quad \gamma=2\alpha-5,\quad 0<\alpha<3,\lb{5.1}\\
&&
b(\cos\theta)=C_{\alpha}
\big((1-\cos\theta)^{-\beta}- (1+\cos\theta)^{-\beta}\big)^2, \quad \beta=(3-\alpha)/2.
\dnumber \lb{5.2}\ees
By elementary calculation we have
\bes&&\big((1-t)^{-\beta}-(1+t)^{-\beta}\big)^2
=\sum_{n=1}^{\infty}a_nt^{2n},\quad |t|<1,\lb{5.3}\\
&&
a_n=4\sum_{i+j=n-1}\fr{\beta(\beta+1)
\cdots(\beta+2i)}{(2i+1)!}\cdot\fr{\beta(\beta+1)\cdots(\beta+2j)}{(2j+1)!}>0,\quad n=1,2,3,...\,.
\dnumber \lb{5.4}\ees
So $b(\cos\theta)$ is completely positive in $\cos^2(\theta)$.
Also we have
\beas&&b(\cos\theta)=\beta^2\fr{\cos^2(\theta)}{\sin^{4\beta}(\theta)}\Big(\int_{-1}^{1}
\big(1+\tau\cos\theta\big)^{\beta-1}
{\rm d}\tau\Big)^2\le 9\fr{\cos^2(\theta)}{\sin^{4\beta}(\theta)}.\eeas
Thus with the function $b_*(\cos\theta)=b^*(\cos\theta)=b(\cos\theta)$ and
$\gm=2\alpha-5$, the kernel
$B(z,\sg)=|z|^{\gm}b(\cos\theta)$ satisfies (A1),(A2),(A4) if
 $5/2\le \alpha<3$, and (A1),(A3),(A4) if
 $3/2\le \alpha<5/2$.
For the case $\alpha=1$, i.e. for the Coulomb potential  $\phi(|x|)=|x|^{-1}$, we have
\be B(z,\sg)={\rm const.}|z|^{-3}\fr{\cos^2(\theta)}{\sin^4(\theta)}\lb{5.4*}\ee
which is essentially equivalent to the Rutherford formula (\ref{Coulomb}) at least for
$0<\sin^2(\theta)<1/2$. In order to satisfy (A3), we have to make an angular cutoff as given in (\ref{Coul}).

 Let
\beas&& G_{t}(|x|)=(2\pi t)^{-\fr{3}{2}}e^{-\fr{|x|^2}{2t}},\quad x\in {\bR},\quad t >0,\\
&& \phi_0(|x|)=\fr{1}{2^{\beta}\Gm(\beta)}\int_{0}^{\infty}G_{t}(|x|) e^{-t/2} t^{\beta-1}{\rm d}t,\quad \beta>0.\eeas
Then
$$\wh{\phi_0}(|\xi|)=(1+|\xi|^2)^{-\beta}.$$
If we choose
$\phi(|x|)=\phi_0(|x|)$ or
$\phi(|x|)=C\dt(x)-\phi_0(|x|)$,
where $C>0$ and $\dt(x)$ is the Dirac delta function at $x=0$, then
$\wh{\phi}(|\xi|)=(1+|\xi|^2)^{-\beta}$ and
$\wh{\phi}(|\xi|)=C-(1+|\xi|^2)^{-\beta}$ respectively and so the corresponding collision kernels $B(z,\sg)$ defined by (\ref{weak-coupling}) are the same:
\beas B(z,\sg)&=&|z|\Big(\big(1+|z|^2\sin^2(\theta/2)\big)^{-\beta}-
\big(1+|z|^2\cos^2(\theta/2)\big)^{-\beta}
\Big)^2\\
&=&\frac{4^{\beta}|z|}{(2+|z|^2)^{2\beta}}\Big(\big(1-a(|z|)\cos\theta\big)^{-\beta}-
\big(1+ a(|z|) \cos\theta\big)^{-\beta}\Big)^2\eeas
with $a(|z|)=\frac{|z|^2}{2+|z|^2}.$
It should be noted that such $B(z,\sg)$ can not be written as a product form
$\Phi(|z|)b(\cos\theta)$.
Using (\ref{5.3}), (\ref{5.4}) we have
\beas&& B(z,\sg)=B(|z|,\cos\theta)=\frac{4^{\beta}|z|}{(2+|z|^2)^{2\beta}}\sum_{n=1}^{\infty}a_n(|z|)\cos^{2n}(\theta),\quad \theta\in(0,\pi),\\
&&
a_n(|z|)=4\sum_{i+j=n-1}\fr{\beta(\beta+1)
\cdots(\beta+2i)}{(2i+1)!}\cdot\fr{\beta(\beta+1)\cdots(\beta+2j)}{(2j+1)!}[a(|z|)]^{2n},\quad n=1,2,3,...\,.\eeas
So $B(z,\sg)$ has the property in Proposition \ref{prop4.1}, i.e.
$B(z,\sg)$ is completely positive in $\cos^2(\theta)$.

On the other hand we have
\beas
B(z,\sg)=\beta^2 |z|^5\cos^2(\theta)\fr{\Big(\int_{0}^{1}
\Big(1+|z|^2\big((1-t)\sin^2(\theta/2)+t\cos^2(\theta/2)\big)\Big)^{\beta-1}
{\rm d}t\Big)^2}{\big\{1+|z|^2+\fr{|z|^4}{4}\sin^2(\theta)\big\}^{2\beta} }.
\eeas
Suppose that $0<\beta<1$. It is easily deduced that
$$|z|^{1-4\beta}\Big( \fr{2|z|^2}{1+|z|^2}\Big)^{2\beta+2}c_{\beta}\cos^2(\theta)
\le B(z,\sg)\le |z|^{1-4\beta}C_{\beta}
\fr{\cos^2(\theta)}{\sin^{4\beta}(\theta) },\quad 0<\beta<1$$
for some constants $0<c_{\beta}<C_{\beta}<\infty$.
Thus, with $\Phi_*(|z|)=(2|z|^2/(1+|z|^2))^{2\beta+2}, b_*(\cos\theta)
=c_{\beta}\cos^2(\theta),  b^*(\cos\theta)=C_{\beta}
\cos^2(\theta)\sin^{-4\beta}(\theta)$, the kernel
$B(z,\sg)$ satisfies (A1),(A2),(A4) for $0\le\gm=1-4\beta<1$ and
(A1),(A3),(A4) for $ -2<\gm =1-4\beta<0$ respectively.
\vskip2mm

{\bf 5.2.} Given any $c> 1,\gm\in{\mR}$ and $\beta>0,\ld>0$, let
\beas&& b(\cos\theta)=c-\cos^2(\theta),\,\,\, \theta\in [0,\pi];\quad
\Phi(|z|)=|z|^{\beta}{\bf 1}_{\{|z|<1\}}+|z|^{\gm}{\bf 1}_{\{|z|\ge 1\}},\,\, z\in{\bR},\\
&& B(z,\sg)=\Phi(|z|)b(\cos\theta),\quad \cos\theta=\la z,\og\ra/|z|,\\
&& h(v)=e^{-\ld|v|^2}\Big(1+\sqrt{2}\sin\big(\fr{\pi}{2}\la \fr{v}{|v|},{\bf e}_1\ra\big)\Big),\quad v\in {\bR}\setminus\{0\}\eeas
 where ${\bf e}_1=(1,0,0)$. Note that the function $b(\cos\theta)
 =c-\cos^2(\theta)$ is positive but not completely positive in $\cos^2(\theta)$.
 We show that for suitably large $\ld>0, \beta>0$,
$$\int_{{\bRRS}}B(z,\sg)h(v)h(v_*)h(v')h(v_*'){\rm d}\sg{\rm d}v{\rm d}v_*< \fr{-5}{\ld^{3/2}\pi^{1/2}}\int_{1}^{\infty}r^{2+\gm} e^{-\ld r^2}{\rm d}r<0 .$$
 To do this we first use
 $e^{-\ld|v|^2}e^{-\ld|v_*|^2}e^{-\ld|v'|^2}e^{-\ld|v_*'|^2}
 =e^{-\ld |v+v_*|^2}e^{-\ld |v-v_*|^2}$ and change of variables
 $(v,v_*)=(\fr{u+z}{2}, \fr{u-z}{2}), z= r\og, u\to r u$ \,to get
\beas&& I:=\int_{{\bRRS}}Bhh_*h'h_*'{\rm d}\sg{\rm d}v{\rm d}v_*
 =\fr{1}{8}\int_{0}^{\infty}\Phi(r)r^5 e^{-\ld r^2}{\rm d}r\int_{\mR^3} e^{-\ld r^2|u|^2}J(u){\rm d}u
 \\
 &&=\fr{1}{8}\int_{0}^{1}r^{5+\beta} e^{-\ld r^2}{\rm d}r\int_{\mR^3} e^{-\ld r^2|u|^2}J(u){\rm d}u
 +\fr{1}{8}\int_{1}^{\infty}r^{5+\gm} e^{-\ld r^2}{\rm d}r\int_{\mR^3} e^{-\ld r^2|u|^2}J(u){\rm d}u
 \eeas
where
\beas J(u)&=&\int_{\bSS}b(\la \og,\sg\ra)
\Big(1+\sqrt{2}\sin\big(\fr{\pi}{2}\la \fr{u+\og}{|u+\og|},{\bf e}_1\ra\big)\Big)
\Big(1+\sqrt{2}\sin\big(\fr{\pi}{2}\la \fr{u-\og}{|u-\og|},{\bf e}_1\ra\big)\Big)
\\
&\times&
\Big(1+\sqrt{2}\sin\big(\fr{\pi}{2}\la \fr{u+\sg}{|u+\sg|},{\bf e}_1\ra\big)\Big)
\Big(1+\sqrt{2}\sin\big(\fr{\pi}{2}\la \fr{u-\sg}{|u-\sg|},{\bf e}_1\ra\big)\Big)
{\rm d}\og {\rm d}\sg.
\eeas
Compute $J(u)$ for $u=0$: using $(1+\sqrt{2}\sin(\theta/2))(1+\sqrt{2}\sin(-\theta/2))
=\cos(\theta)$ we have
\beas J(0)
&=&\int_{\bSS}\big(c- \la \og,\sg\ra^2\big)\cos\big(\pi\la \og,{\bf e}_1\ra)\cos\big(\pi\la \sg,{\bf e}_1\ra)
{\rm d}\og {\rm d}\sg
\\
&=&c\Big(\int_{\bS}\cos(\pi\la\og,e_1\ra){\rm d}\og\Big)^2
-\int_{\bSS}\la \og,\sg\ra^2\cos\big(\pi\la \og,{\bf e}_1\ra)\cos\big(\pi\la \sg,{\bf e}_1\ra)
{\rm d}\og {\rm d}\sg\\
&=&c\Big(2\pi\int_{-1}^1\cos(\pi t){\rm d}t\Big)^2-\sum_{i=1}^3\sum_{j=1}^3\Big(\int_{\bS}\og_i\og_j\cos(\pi\og_1){\rm d}\og\Big)^2
\quad (\og=(\og_1, \og_2,\og_3))\quad \\
&=&0 -\sum_{i=1}^3\Big(\int_{\bS}\og_i^2\cos(\pi\og_1){\rm d}\og\Big)^2
=-\fr{96}{\pi^2}
\eeas
The function $J(u)$ is continuous on
$|u|\le 1/2$. So there is $0<\dt\le 1/2$ such that
$J(u)\le -\fr{90}{\pi^2}$ for all $|u|\le \dt$.
Also we have $|J(u)|\le C_0:=c(1+\sqrt{2})^4(4\pi)^2$
for all $u\in {\bR}$. Thus
\beas&&\int_{\mR^3} e^{-\ld r^2|u|^2}J(u){\rm d}u
\le -\fr{90}{\pi^2}\int_{|u|\le \dt} e^{-\ld r^2|u|^2}{\rm d}u
+C_0\int_{|u|>\dt} e^{-\ld r^2|u|^2}{\rm d}u
\\
&&=\ld ^{-3/2}r^{-3}\Big( -\fr{90}{\pi^{1/2}}
+\big(\fr{90}{\pi^2}+ C_0\big)\int_{|u|>\dt\sqrt{\ld} r} e^{-|x|^2}{\rm d}x
\Big),\quad r\ge 1.
\eeas
Choose $\ld>0$ large enough  such that
$$
(\fr{90}{\pi^2}+ C_0)\int_{|u|>\dt\sqrt{\ld}} e^{-|x|^2}{\rm d}x
\le  \fr{10}{\pi^{1/2}}.$$
Then for this $\ld$
\beas&&\int_{\mR^3} e^{-\ld r^2|u|^2}J(u){\rm d}u
\le -\ld ^{-3/2}r^{-3}\fr{80}{\pi^{1/2}}\qquad \forall\, r\ge 1,\\
&&
\fr{1}{8}\int_{1}^{\infty}\{\cdots\}\le -\ld^{-3/2}\fr{10}{\pi^{1/2}}\int_{1}^{\infty}r^{2+\gm} e^{-\ld r^2}{\rm d}r.
\eeas
On the other hand,
$$\fr{1}{8}\int_{0}^{1}\{\cdots\}\le \ld^{-3/2}C_0 \fr{\pi^{3/2}}{8}\int_{0}^{1}r^{2+\beta} e^{-\ld r^2}{\rm d}r
\le  \ld^{-3/2}C_0 \fr{\pi^{3/2}}{8}\cdot\fr{1}{3+\beta}.$$ Choose
$\beta>0$ large enough such that $\fr{1}{8}\int_{0}^{1}\{\cdots\}<\ld^{-3/2}\fr{5}{\pi^{1/2}}\int_{1}^{\infty}r^{2+\gm} e^{-\ld r^2}{\rm d}r$.
 Then for the chosen $\ld,\beta$ we
 obtain
 $$I=\fr{1}{8}\int_{0}^{1}\{\cdots\}+\fr{1}{8}\int_{1}^{\infty}\{\cdots\}<-\ld^{-3/2}\fr{5}{\pi^{1/2}}\int_{1}^{\infty}r^{2+\gm} e^{-\ld r^2}{\rm d}r. $$
\vskip6mm

{\bf Acknowledgment}. This work was supported by NSF of China under Grant No.11771236.\\

\noindent {\bf Data Availability}\, Data sharing is not applicable to this article as no datasets were generated or analysed during the current study.
\vskip2mm

\noindent {\bf Declarations}

\noindent {\bf Conflict of Interest}\, The authors declare no conflict of interest.
\\


\begin{thebibliography}{99}

\bibitem{Alex} Alexandre, R.: On some related non homogeneous 3D Boltzmann models in the non cutoff case. J. Math. Kyoto Univ.{\bf 40}, no. 3, 493-524 (2000).

\bibitem{ABDL} Alonso,R.; Bagland,V.; Desvillettes,L.; Lods,B.: About the Landau-Fermi-Dirac equation with moderately soft potentials. Arch. Ration. Mech. Anal. {\bf 244}, no.3, 779-875 (2022).

\bibitem{weak-coupling} Benedetto, D., Pulvirenti, M., Castella, F., Esposito, R.: On
the weak-coupling limit for bosons and fermions. Math. Models Methods
Appl. Sci. {\bf 15}, 1811-1843  (2005).

\bibitem{Bohm} Bohm, D.: Quantum Theory, Prentice-Hall, New York, 1951.

\bibitem{BD} Bouchut, F.; Desvillettes, L.: A proof of the smoothing properties of the positive part of Boltzmann's kernel. Rev. Mat. Iberoamericana {\bf 14}, no. 1, 47-61 (1998).

\bibitem{CCL2009} Carlen, E.A., Carvalho, M.C., Lu, X.: On strong convergence to equilibrium for the Boltzmann
    equation with soft potentials. J. Stat. Phys. {\bf 135}, 681-736 (2009)

\bibitem {Chapman and Cowling} Chapman, S., Cowling, T.G.: {\it The Mathematical
Theory of Non-Uniform Gases}. Third Edition (Cambridge University
Press, 1970).

\bibitem{D} Desvillettes, L.: Some applications of the method of moments for the homogeneous Boltzmann and Kac equations. Arch. Rational Mech. Anal. {\bf 123}, no.4, 387-404 (1993).

\bibitem{Do} Dolbeault, J.: Kinetic models and quantum effects: a modified Boltzmann equation for Fermi-Dirac particles. Arch. Rational Mech. Anal. {\bf 127}, no. 2, 101-131 (1994).


\bibitem{DS} Dunford, N.; Schwartz,J.T.: {\it Linear Operators I: General Theory} (Interscience, New York,
1958).


\bibitem{ESY} Erd\"{o}s, L., Salmhofer, M., Yau, H.-T.: On the
quantum Boltzmann equation. J. Stat. Phys. {\bf 116}, 367-380  (2004).


\bibitem{EMV0} Escobedo, M., Mischler, S., Valle, M.A.:
Homogeneous Boltzmann equation in quantum relativistic kinetic theory.
Electronic Journal of Differential Equations, Monograph, 4.
Southwest Texas State University, San Marcos, TX, 2003. 85 pp.


\bibitem{Lu1998} Lu, X.: A direct method for the regularity of the gain term in the Boltzmann equation. J. Math. Anal. Appl. {\bf228}, no. 2, 409-435 (1998).

\bibitem{Lu2001} Lu, X.: On spatially homogeneous solutions of a modified Boltzmann equation for Fermi-Dirac particles. J. Statist. Phys. {\bf 105}, no. 1-2, 353-388 (2001).

 \bibitem{Lu2008} Lu, X.: On the Boltzmann equation for Fermi-Dirac particles with very soft potentials: global existence of weak solutions.
     J.Differential Equations {\bf 245}, no.7, 1705-1761 (2008).

\bibitem{LW2003} Lu, X., Wennberg, B.: On stability and strong convergence for the spatially homogeneous Boltzmann equation for Fermi-Dirac particles. Arch. Ration. Mech. Anal.
     {\bf 168}, no.1, 1-34 (2003).

\bibitem{LS} Lukkarinen, J., Spohn, H.:
Not to normal order--notes on the kinetic limit for weakly interacting quantum fluids.
 J. Stat. Phys. {\bf 134}, 1133-1172  (2009).

\bibitem {Nordheim} Nordheim, L.W.: On the kinetic methods in the new statistics
and its applications in the electron theory of conductivity.
Proc. Roy. Soc. London Ser. A  {\bf 119}, 689-698  (1928)

\bibitem {OW} Ouyang, Z.; Wu, L.: On the quantum Boltzmann equation near Maxwellian and vacuum. J. Differential Equations {\bf 316}, 471-551  (2022).

\bibitem{Lions1} Lions, P.-L.: Compactness in Boltzmann's equation via Fourier integral operators and applications. I, II. J. Math. Kyoto Univ. {\bf 34}, no. 2, 391-427, 429-461  (1994).

\bibitem{Lions2} Lions, P.-L.: Compactness in Boltzmann's equation via Fourier integral operators and applications. III. J. Math. Kyoto Univ. 34, no. 3, 539-584 (1994).

\bibitem {Stein} Stein, Elias M.: {\it Singular integrals and differentiability properties of functions}. Princeton Mathematical Series, No.30 Princeton University Press, Princeton, N.J. 1970 xiv+290 pp.

\bibitem {Uehling and Uhlenbeck}Uehling, E.A., Uhlenbeck, G.E.:
Transport phenomena in Einstein-Bose and Fermi-Dirac gases, I, Phys. Rev.
{\bf 43}, 552-561  (1933).

\bibitem{Villani2} Villani, C.: A review of mathematical topics in collisional kinetic theory.  Handbook of mathematical fluid dynamics, Vol.I, 71-305, North-Holland, Amsterdam, 2002.

\bibitem{Wen1997}  Wennberg, B.: The geometry of binary collisions and generalized Radon transforms. Arch. Rational Mech. Anal. {\bf 139}, no.3, 291-302 (1997).





\end{thebibliography}
\end{document}